\documentclass[a4paper,twoside,10pt]{article}

\usepackage[a4paper,left=3cm,right=3cm, top=3cm, bottom=3cm]{geometry}
\usepackage[latin1]{inputenc}
\usepackage{mathrsfs}
\usepackage{graphicx}
\usepackage{epstopdf}
\usepackage{amsmath}
\usepackage{amsthm}
\usepackage{amssymb}
\usepackage{hyperref}
\usepackage{stmaryrd}
\usepackage{subfigure}
\usepackage{float}
\usepackage{bigints}
\usepackage{cite}
\usepackage{color}
\usepackage[abs]{overpic}
\usepackage[font=footnotesize,labelfont=bf]{caption}
\usepackage{cases}
\usepackage{tikz}
\usepackage{algorithmic}
\usepackage{rotating}
\usepackage{blkarray}
\usetikzlibrary{arrows}
\usetikzlibrary{quotes,angles}
\usepackage[author={Lorenzo}]{pdfcomment}
\usepackage{soul,xcolor}

\restylefloat{table}
\theoremstyle{plain}

\theoremstyle{definition}

\theoremstyle{remark}

\newcommand{\boldalpha}{\boldsymbol{\alpha}} 
\newcommand{\boldbeta}{\boldsymbol{\beta}} 
\newcommand{\boldgamma}{\boldsymbol{\gamma}} 
\newcommand{\p}{p} 
\renewcommand{\u}{u} 
\renewcommand{\v}{v} 
\newcommand{\f}{f} 
\newcommand{\fn}{\f_n} 
\newcommand{\V}{V} 
\renewcommand{\a}{a} 
\newcommand{\aE}{\a^\E} 
\newcommand{\an}{\a_\p} 
\newcommand{\anE}{\a_\p^\E} 
\newcommand{\taun}{\mathcal T_n} 
\newcommand{\E}{K} 
\newcommand{\En}{\E_i} 
\newcommand{\NE}{N_\E} 
\newcommand{\e}{e} 
\newcommand{\Pizpmd}{\Pi^0_{\p-2}} 
\newcommand{\Pinablap}{\Pi^{\nabla}_\p} 
\newcommand{\un}{\u_\p} 
\newcommand{\unn}{\u_\p} 
\newcommand{\vnn}{\v_\p} 
\newcommand{\Vn}{\V_\p} 
\newcommand{\VnE}{\Vn(\E)} 
\newcommand{\VnuE}{\V_{\p,1}(\E)} 
\newcommand{\VndE}{\V_{\p,2}(\E)} 
\newcommand{\q}{q} 
\newcommand{\qpmd}{\q_{\p-2}} 
\newcommand{\qp}{\q_\p} 
\newcommand{\qalpha}{\q_{\boldalpha}} 

\newcommand{\qbeta}{\q_{\boldbeta}} 
\newcommand{\qgamma}{\q_{\boldgamma}} 
\newcommand{\qalphau}{\q_{\boldalpha}^1} 
\newcommand{\qbetau}{\q_{\boldbeta}^1} 
\newcommand{\qalphad}{\q_{\boldalpha}^2} 
\newcommand{\qbetad}{\q_{\boldbeta}^2} 
\newcommand{\qgammad}{\q_{\boldgamma}^2} 
\newcommand{\qalphat}{\q_{\boldalpha}^3} 
\newcommand{\azpmd}{_{\vert \boldalpha \vert=0}^{\p-2}} 
\newcommand{\aupmd}{_{\vert \boldalpha \vert=1}^{\p-2}} 
\newcommand{\azp}{_{\vert \boldalpha \vert=0}^{\p}} 
\newcommand{\dof}{\text{dof}} 
\newcommand{\qalphabasis}{\{\qalpha\}} 
\newcommand{\qalphabasisi}{\{\qalpha^i\}} 
\newcommand{\qalphabasisu}{\{\qalpha^1\}} 
\newcommand{\qalphabasisd}{\{\qalpha^2\}} 
\newcommand{\qalphabasist}{\{\qalpha^3\}} 
\newcommand{\n}{\mathbf n} 
\newcommand{\SE}{S^\E} 
\newcommand{\SEmat}{\mathbf {\SE}} 
\newcommand{\SEmatthree}{\mathbf {S^\E_3}} 
\renewcommand{\c}{c} 
\newcommand{\upi}{\u_\pi} 
\newcommand{\uI}{\u_I} 
\newcommand{\Fn}{\mathcal F_n} 
\newcommand{\h}{h} 
\newcommand{\hE}{\h_\E} 
\newcommand{\NdofE}{N_{\dof}^\E}
\newcommand{\Ndofbndr}{N_{\dof}^{\text{bndr}}} 
\newcommand{\xbold}{\mathbf x} 
\newcommand{\xEbold}{\mathbf{x_\E}} 
\newcommand{\xE}{x_\E} 
\newcommand{\yE}{y_\E} 
\newcommand{\GS}{\mathbf{GS}} 
\newcommand{\G}{\mathbf G} 
\newcommand{\Gtilde}{\widetilde{\G}} 
\newcommand{\A}{\mathbf A} 
\newcommand{\B}{\mathbf B} 
\newcommand{\D}{\mathbf D} 
\renewcommand{\H}{\mathbf H} 
\newcommand{\C}{\mathbf C} 
\renewcommand{\L}{\mathbf L} 
\newcommand{\Gbar}{\overline{\G}} 
\newcommand{\Gtildebar}{\overline {\Gtilde}} 
\newcommand{\Bbar}{\overline \B} 
\newcommand{\Dbar}{\overline \D} 
\newcommand{\Cbar}{\overline \C} 
\newcommand{\Lbar}{\overline \L} 
\newcommand{\Fbar}{\overline{\mathbf F}} 
\newcommand{\npmd}{n_{\p-2}} 
\newcommand{\np}{n_\p} 
\newcommand{\Id}{\mathbf {Id}}  
\newcommand{\K}{\mathbf K} 
\newcommand{\Pinablastar}{\mathbf {\Pi^{\nabla}_*}} 
\newcommand{\Pinablamat}{\mathbf {\Pi^{\nabla}}} 

\begin{tiny}
\author{
\normalsize{
L. Mascotto
\thanks{Faculty of Mathematics, University of Vienna, 1090 Vienna, Austria. E-mail: {\tt lorenzo.mascotto@univie.ac.at}}
}}
\end{tiny}

\date{}

\title{\textbf{\normalsize{Ill-conditioning in the Virtual Element Method: stabilizations and bases}}}

\begin{document}

\maketitle

\begin{abstract}
In this paper we investigate the behavior of the condition number of the stiffness matrix resulting from the approximation of a 2D Poisson problem by means of the Virtual Element Method.
It turns out that ill-conditioning appears when considering high-order methods or in presence of ``bad-shaped'' (for instance nonuniformly star-shaped, with small edges\dots) sequences of polygons.
We show that in order to improve such condition number one can modify the definition of the internal moments by
choosing proper polynomial functions that are not the standard monomials.
We also give numerical evidence that, at least for a 2D problem, 
standard choices for the stabilization give similar results in terms of condition number.
\end{abstract}

\section {Introduction} \label{section introduction}
The interest in numerical methods for the approximation of Partial Differential Equations (PDEs in short) based on polytopic grids has grown in the last decade, due to the high flexibility that polygonal/polyhedral meshes allow.
Among the other methods, we here recall only a short list including:
Mimetic Finite Differences \cite{BLM_MFD, BLS_MFD}, Discontinuous Galerkin-Finite Element Method (DG-FEM) \cite{antonietti2016reviewDG, cangianigeorgoulishouston_hpDGFEM_polygon},
Hybridizable and Hybrid High-Order Methods \cite{dipietroErn_hho, cockburn2008superconvergent}, Weak Galerkin Method \cite{wang2013weak},
BEM-based FEM \cite{Weisser_basic} and Polygonal FEM \cite{SukumarTabarraeipolygonalintroduction}.

An alternative approach is offered by the Virtual Element Method (VEM in short), recently introduced in \cite{VEMvolley}.
VEM are a generalization of the Finite Element Method (FEM in short) enabling the employment of polytopal meshes and the possibility of building high-order methods.
In addition to polynomials, local VE spaces consist of other functions instrumental for constructing global $H^1$ conforming approximation space; such functions are tipically the solution of local PDEs and therefore are not known explicitely in a closed-form.
For this reason, VEM can be considered, for all practical purposes, Trefftz methods.

In VEM, the bilinear forms and the right-hand sides are not computed exactly due to the fact that the functions in VE spaces are not fully explicit.
The VEM gospel states that such bilinear forms and right-hand sides are then approximated by means of discrete counterparts that are exactly computable only through the degrees of freedom and which  scale like the continuous ones.

In a few years, thanks to high theoretical and practical flexibility of VEM, the size of associated literature has rapidly blown up. Among the other references, we recall the following:
the $\p$ and $\h\p$ version of the method \cite{hpVEMbasic, preprint_hpVEMcorner, pVEMmultigrid, preprint_HarmonicVEM, fetishVEM3D}, parabolic problems \cite{parabolicVEM},
Cahn-Hilliard, Stokes, Navier-Stokes and Helmoltz equations \cite{absv_VEM_cahnhilliard, BLV_StokesVEMdivergencefree, Helmholtz-VEM, VEM_NavierStokes,streamvirtualelementformulationstokesproblempolygonalmeshes},
linear and nonlinear elasticity problems \cite{BLM_VEMsmalldeformation, VEMelasticity, Paulino-VEM}, general elliptic problems \cite{BBMR_generalsecondorder}, PDEs on surfaces \cite{VEM_LaplaceBeltrami},
Domain Decomposition \cite{VEM_DD_basic}, 
application to discrete fracture networks \cite{Benedetto-VEM-3}, serendipity VEM \cite{serendipityVEM}, VEM on surfaces \cite{VEM_LaplaceBeltrami}. The implementation of the method is described in \cite{hitchhikersguideVEM},
whereas the basic principles of the 3D version of the method are the topic of \cite{equivalentprojectorsforVEM, preprint_VEM3Dbasic}.

It was observed that the VEM stiffness matrix, similarly as FEM, can become ill-conditioned in various situations.
This is the case for instance of the $\p$ and the $\h\p$ version of VEM, as discussed in \cite{hpVEMbasic, pVEMmultigrid}, and in presence of ``bad-shaped'' (i.e. nonuniformly star-shaped, with small edges\dots) polygons, as discussed \cite{BerroneBorio_orthogonalVEM}.

Among the possible reasons of this ill-conditioning we highlight two of them. The first one is related to the fact that in the VEM framework one does not employ the exact bilinear form but an approximated one;
the choice of the discrete bilinear form, and, in particular, of one of its two components, namely the stabilization of the method, may have an impact on the 3D version of VEM as observed in \cite{preprint_VEM3Dbasic}.
The second one is the choice of the basis. This is also the case for FEM, where the choice of the basis has an important role on the ill-conditioning of the system, see \cite{SchwabpandhpFEM, dubinerspectraltriangles, AdjeridAiffaFlahertyhierarchicalfiniteelementbasis}
and the references therein.

The aim of the present paper is to discuss various choices for both the discrete bilinear forms and the VE bases and check numerically
that particular choices can cure the ill-conditioning which arises in high-order (or in presence of bad-shaped polygons) VEM.

In particular, we show that, while the choice of the stabilization has not a deep impact on the ill-conditioning (at least for the 2D case, which is the focus of this paper),
a proper choice of the basis can actually improve the condition number of the stiffness matrix.
However, it is worth to mention that in ``standard'' situations, i.e. low-to-moderate order VEM and VEM applied to shape-regular decompositions, the standard VEM (e.g. the one described in \cite{VEMvolley})
is preferable, since the implementation of that version of the method turns out to be much simpler than those we are going to present in this article.

The outline of the paper is the following. We firstly discuss the model problem and its VEM approximation in Section \ref{section VEM},
while in Appendix \ref{section hitchhikers} we give a hint on the implementation details of the method using one of the new VEM basis.
In Section \ref{section numerical results}, we present a number of numerical experiments comparing the behaviour of the method when changing various stabilizations and VE bases.

In the remainder of the paper, we adopt the standard notation for Sobolev spaces, see e.g. \cite{evansPDE, adamsfournier}.
Given $\omega \subset \mathbb R^2$, we denote by $H^\ell (\omega)$, $\ell \in \mathbb N$, the Sobolev space of order $\ell$ over $\omega$;
in the case $\ell = 0$, we set $H^0(\omega) = L^2(\omega)$, where $L^2(\omega)$ is the Lebesgue space of square integrable functions over $\omega$.
By $H^1_0(\Omega)$, we mean the Sobolev space $H^1$ of functions with zero trace. The Sobolev (semi)norms and inner products read, respectively:
\begin{equation} \label{Sobolev stuff}
\vert \cdot \vert_{\ell, \omega},\quad \quad \quad \quad \Vert \cdot \Vert_{\ell, \omega},\quad \quad \quad  \quad(\cdot, \cdot)_{0,\omega}.
\end{equation}
Further, given a Lipschitz domain $\omega$ with boundary $\Gamma$, we set $\n$ the normal versor associated with $\Gamma$ and $\partial _\n v = \nabla \v \cdot \n$ the normal derivative of a sufficiently regular function $v$.
By $\mathbb P_\ell (\omega)$, $\ell \in \mathbb N$, we denote the space of polynomials of degree $\ell$ over $\omega$.
Finally, given two positive quantities $a$ and $b$, possibly depending on the discretization parameters $\h$ and $\p$,
we write $a \lesssim b$ if there exists a positive and parameter-independent constant $c$ such that $a \le c \, b$.
Moreover, we write $a \approx b$ meaning that $a \lesssim b$ and $b \lesssim a$ simultaneously.

\section {The Virtual Element Method: definition, choice of the stabilization and of the basis} \label{section VEM}
Given $\Omega \subset \mathbb R^2$ a polygonal domain and $\f \in L^2(\Omega)$, we consider the following 2D Poisson problem:
\begin{equation} \label{Poisson problem}
\begin{cases}
\text{find } \u \in \V \text{ such that}\\
\a(\u,\v) = (\f,\v)_{0,\Omega} \quad \forall \v \in \V\\
\end{cases},
\quad \text{where } \V := H^1_0(\Omega),\quad \text{and} \quad \a(\cdot, \cdot) :=  (\nabla \cdot, \nabla \cdot)_{0,\Omega}.
\end{equation}
The outline of the present section is the following. In Section \ref{subsection a family of VEM}, we introduce a family of VEMs for the approximation of problem \eqref{Poisson problem};
here, both the stabilization, typical of VEM, and the choice of the basis of local VE space are kept at a very general level.
Explicit choices for the stabilization are investigated in Section \ref{subsection choices for the stabilization}, whereas explicit choices for the local basis elements are the topic of Section \ref{subsection choices for the basis}.
Finally, in Section \ref{subsection stab and bases effects on the method}, we highlight the influence of the choices of the stabilizations and of the local basis elements on the performances of the method.

\subsection{A family of VEM} \label{subsection a family of VEM}
Given $\Omega$ the computational domain of problem \eqref{Poisson problem}, we consider a family $\{\taun\}_{n \in \mathbb N}$ of conforming polygonal decomposition of $\Omega$,
where by conforming we mean that, given an edge in the skeleton of the decomposition which does not lie on $\partial \Omega$, then it is an edge of \emph{exactly} two polygons.
We also fix $\p \in \mathbb N$, which will represent the degree of accuracy of the method.
We associate to each $\E \in \taun$ its diameter $\hE$ and to decomposition $\taun$ its mesh size function $\h = \max_{\E \in \taun} \hE$.

Standard regularity assumptions on decomposition $\taun$ that are usually required in VEM literature read:
\begin{itemize}
\item[(\textbf{D1})] For every $\E \in \taun$, $\E$ is star-shaped (see \cite{BrennerScott}) with respect to a ball of radius greater or equal than $\gamma \hE$, $\gamma$ being a positive constant independent of the family of decompositions.
\item[(\textbf{D2})] For every $\E \in \taun$ and for every $\e$ edge of $\E$, the length of $\e$ is greater or equal than $\gamma \hE$, $\gamma$ being the same constant introduced in assumption (\textbf{D1}).
\end{itemize}
We point out that such assumptions can be relaxed, see \cite{beiraolovadinarusso_stabilityVEM}.

We now define the local VE space on polygon $\E \in \taun$ following the standard definition given e.g. in \cite{VEMvolley}. Having set the space of piecewise continuous polynomials of degree $\p$ over $\partial \E$:
\begin{equation} \label{boundary space}
\mathbb B_\p(\partial \E) := \left\{  \vnn \in \mathcal C^0(\partial \E) \mid \vnn |_\e \in \mathbb P_\p(\e) \; \forall \e \text{ edge of } \E  \right\},
\end{equation}
we introduce the local VE space as:
\begin{equation} \label{local VE space}
\VnE := \left\{ \vnn \in H^1(\E) \mid \vnn |_{\partial \E} \in  \mathbb B_\p(\partial \E),\; \Delta \vnn \in \mathbb P_{\p-2}(\E)   \right\}.
\end{equation}
Space $\VnE$, which contains $\mathbb P_\p(\E)$, contains also other functions that in general are not known pointwise (hence, the name \emph{virtual}),
but which are added to polynomials in order to guarantee the possibility of building a $H^1$ conforming method over $\taun$.

We associate to space $\VnE$ the following set of linear functionals. Given $\vnn \in \VnE$:
\begin{itemize}
\item the point values of $\vnn$ at the vertices of $\E$;
\item for every $\e$ edge of $\E$, the point values of $\vnn$ at the $\p-1$ internal Gau\ss-Lobatto nodes of $\e$;
\item the (scaled) internal moments:
\begin{equation} \label{internal moments}
\frac{1}{\vert \E \vert} \int_\E \vnn \, \q_\alpha,
\end{equation}
where $\{ \q_\alpha \} _{\alpha=1}^{\dim(\mathbb P_{\p-2}(\E))}$ is \emph{any} basis of $\mathbb P_{\p-2}(\E)$.
\end{itemize}
It was proven in \cite{VEMvolley} that this is a set of unisolvent degrees of freedom.
Let $\NdofE$ be the number of such degrees of freedom, i.e. the dimension of space $\VnE$.
Henceforth, we denote by $\dof_i$ the $i$-th dof of space $\VnE$ and by:
\begin{equation} \label{canonical basis}
\{\varphi_i\} _{i=1}^{\NdofE}
\end{equation}
the local canonical basis of space $\VnE$, which is dual to the set of degrees of freedom $\{\dof_i\}_{i=1}^{\dim(\VnE)}$, i.e. the set of functions such that:
\begin{equation} \label{dual basis}
\dof_i(\varphi_j) = \delta_{i,j} = \begin{cases}
1 & \text{if } i = j\\
0 & \text{otherwise}\\
\end{cases}\quad \quad \forall i,\, j = 1, \dots, \NdofE.
\end{equation}
Functions in the canonical basis associated with internal moments \eqref{internal moments} are bubbles on element $\E$ since they vanish on the boundary.

It is fundamental to observe that the definition of space $\VnE$ is completely independent of the choice of polynomial basis $\{ \q_\alpha \} _{\alpha=1}^{\dim(\mathbb P_{\p-2}(\E))}$,
which is, so far, \emph{only} instrumental for the definition of the internal moments \eqref{internal moments}.
Nonetheless, such a choice plays a crucial role in the behaviour of the condition number of the stiffness matrix of the method as discussed and shown in Section \ref{section numerical results}.
In Subsection \ref{subsection choices for the basis}, we present many choices of basis $\qalphabasis \azpmd$.

As already stressed, not all the functions in space $\VnE$ are known explicitly; nevertheless, it is possible to compute local $H^1$ and $L^2$ projections as shown in \cite{VEMvolley, hitchhikersguideVEM} via the degrees of freedom above described.
More precisely, it is possible to compute the following two operators.

The first one is the $L^2$ projector $\Pizpmd : \VnE \rightarrow \mathbb P_{\p-2}(\E)$ defined, for all $\E \in \taun$, by:
\begin{equation} \label{L2 orthogonality}
(\qpmd, \Pizpmd \vnn - \vnn)_{0,\E} = 0 \quad \forall \qpmd \in \mathbb P_{\p-2}(\E),\; \forall \vnn \in \VnE,
\end{equation}
which is clearly computable via internal moments \eqref{internal moments}.

The second one is a $H^1$ projector $\Pinablap : \VnE \rightarrow \mathbb P_\p(\E)$ defined by:
\begin{equation} \label{H1 orthogonality}
\begin{cases}
\aE(\qp, \Pinablap \vnn - \vnn) = 0\\
P_0 (\Pinablap \vnn - \vnn) = 0\\
\end{cases}
\quad \forall \qp \in \mathbb P_\p(\E),\; \forall \vnn \in \VnE,
\end{equation}
where $\aE(\cdot, \cdot) = (\nabla \cdot, \nabla \cdot)_{0,\E}$ and where $P$ is an operator which fixes the constant part of the energy projector $\Pinablap$ defined as:
\begin{equation} \label{fixing constants}
P_0(\vnn) = \begin{cases}
\frac{1}{\NE} \sum_{j=1}^{\NE} \vnn (\nu_j) & \text{if } \p=1\\
\int_{\E} \vnn & \text{otherwise}\\
\end{cases}\quad \forall \vnn \in \VnE,
\end{equation}
where $\{\nu_j\}_{j=1}^{\NE}$ is the set of vertices of polygon $\E$.
Operator $\Pinablap$ is computable by means of the degrees of freedom, noting that:
\[
\aE(\qp, \vnn) = - \int_\E \Delta \qp \; \vnn + \int_{\partial \E} \partial _\n \qp \; \vnn \quad \forall \qp \in \mathbb P_\p (\E),\; \forall \vnn \in \VnE
\]
and recalling that $\vnn \in \VnE$.

At this point, we turn our attention to the computation of the local stiffness matrix and of the local right-hand side. In both cases, we can not use their continuous counterparts since we do not know pointwise functions in space $\VnE$.
Therefore, we follow the VEM gospel and we note that Pythagoras Theorem for Hilbert spaces asserts:
\begin{equation} \label{Pythagoras theorem}
\aE(\unn, \vnn) = \aE (\Pinablap \unn, \Pinablap \vnn) + \aE( (I-\Pinablap) \unn, (I-\Pinablap) \vnn)\quad \forall \unn,\, \vnn \in \VnE.
\end{equation}
If we split the local space $\VnE$ into a polynomial and a ``pure virtual'' part:
\begin{equation} \label{splitting space VnE}
\begin{split}
\VnE 	& = \left\{ \vnn \in \VnE \mid \vnn \in \mathbb P_\p(\E)  \right\} \oplus \left\{ \vnn \in \VnE \mid \Pinablap \vnn = 0 \right\} \\
	& = \mathbb P_\p(\E) \oplus \ker(\Pinablap) =: \VnuE \oplus \VndE,
\end{split}
\end{equation}
then the first term in the right-hand side of \eqref{Pythagoras theorem} is identically zero on $\VndE$, while the second one annihilates on $\VnuE$.

We point out that, given $\NE$ the number of vertices of $\E$, one has:
\begin{equation} \label{dimension of subspaces}
\begin{split}
& \dim (\VnuE) = \dim(\mathbb P_\p(\E)) = \frac{(\p+1)(\p+2)}{2},\\
\medskip
& 
\begin{split}
\dim (\VndE) &= \dim(\VnE) - \dim(\VnuE) \\
		&= \NE \cdot \p + \frac{(\p-1)\p}{2} - \frac{(\p+1)(\p+2)}{2} = (\NE-2) \p - 1.\\
\end{split}\\
\end{split}
\end{equation}
This entails that the actual ``pure virtual'' part of space $\VnE$, i.e. $\ker(\Pinablap)$, is asymptotically smaller than its polynomial counterpart,
if the number of vertices of $\E$ remains uniformly bounded. More precisely, $\dim(\VndE) \approx \p$ whereas $\dim(\VnuE) \approx \p^2$.

We emphasize that the first term in the right-hand side of \eqref{Pythagoras theorem} is explicitly computable but the second one is not. For this reason, one substitutes:
\[
\aE((I-\Pinablap) \unn , (I-\Pinablap) \vnn) \quad \Longrightarrow \quad \SE((I-\Pinablap) \unn, (I-\Pinablap) \vnn),
\]
where $\SE(\cdot, \cdot)$ is \emph{any} bilinear form mimicking the continuous one, i.e. $\aE(\cdot, \cdot)$ on $\ker(\Pinablap) \times \ker (\Pinablap)$.
In order to have a well-posed method, we demand the following stability assumption on the bilinear form $\SE$:
\begin{equation} \label{stability S}
\c_*(\p) \vert \vnn \vert^2_{1,\E} \le \SE(\vnn, \vnn) \le \c^*(\p) \vert \vnn \vert^2_{1,\E}\quad \forall \vnn \in \ker(\Pinablap),
\end{equation}
where $\c_*(\p)$ and $\c^*(\p)$ are two positive constants possibly depending on $\p$.
Various possible stabilizations $\SE$ are presented in Section \ref{subsection choices for the stabilization}.

By defining the local discrete bilinear form as:
\begin{equation} \label{local discrete bilinear form}
\anE(\unn, \vnn) = \aE(\unn, \vnn) + \SE((I-\Pinablap) \unn, (I-\Pinablap) \vnn) \quad \forall \unn,\; \vnn \in \VnE,
\end{equation}
one can prove 
the following properties for $\anE(\cdot, \cdot)$:
\begin{itemize}
\item[(\textbf{P1})] \textbf{$\p$-consistency}: for every $\qp \in \mathbb P_{\p}(\E)$ and for every $\vnn \in \VnE$, the local bilinear form $\anE$ satisfies:
\begin{equation} \label{consistency}
\aE(\qp, \vnn) = \anE (\qp, \vnn);
\end{equation}
\item[(\textbf{P2})] \textbf{stability}: for every $\vnn \in \VnE$, the local bilinear form $\anE$ satisfies:
\begin{equation} \label{stability}
\alpha_*(\p) \vert \vnn \vert^2 \le \anE(\vnn, \vnn) \le \alpha^*(\p) \vert \vnn \vert^2,
\end{equation}
where $\alpha_*(\p) = \min (1,\c_*(\p))$ and $\alpha^*(\p) = \max(1,\c^*(\p))$.
\end{itemize}
Regarding the local discrete right-hand side, we set:
\begin{equation} \label{local discrete rhs}
\langle \fn, \vnn \rangle_\E :=
\begin{cases}
\int_{\E} \left( \f \int_{\partial \E} \vnn  \right)  & \text{if } \p= 1\\
\int_\E \f \, \Pizpmd \vnn & \text{otherwise}\\
\end{cases}
\quad \forall \vnn \in \VnE.
\end{equation}
At this point, we are able to define the global VE space, which is obtained by a standard dof coupling of the local spaces:
\begin{equation} \label{global VE space}
\Vn := \left\{ \vnn \in \V \cap \mathcal C^0(\overline \Omega) \mid \vnn |_\E \in \VnE \; \forall \E \in \taun   \right\}
\end{equation}
and the discrete global bilinear form and right-hand side:
\begin{equation} \label{global bilinear form and rhs}
\an(\unn, \vnn) = \sum_{\E \in \taun} \anE(\unn,\vnn),\quad \quad \langle \fn, \vnn \rangle = \sum_{\E \in \taun}  \langle \fn, \vnn \rangle_{\E}\quad  \quad \forall \unn,\,\vnn \in \Vn.
\end{equation}
We now define a family of VEMs depending on the choice of the local stabilization $\SE$ defined in \eqref{stability S}:
\begin{equation} \label{family of VEM}
\begin{cases}
\text{find } \un \in \Vn \text{ such that }\\
\an(\un, \vnn) = \langle \fn, \vnn \rangle\quad \forall \vnn \in \Vn
\end{cases}.
\end{equation}
Property (\textbf{P1}) guarantees that the VEM passes the patch test for piecewise-polynomial of degree $\p$ solutions, i.e. if the solution of problem \eqref{Poisson problem} is a polynomial of degree $\p$
and the order of the method is also $\p$, then the method returns as an output such a polynomial. For this a reason, $\p$ is also called the \emph{degree of accuracy} of the method.

On the other hand, property (\textbf{P2}) implies the coercivity and the continuity of the discrete bilinear form and thus the well-posedness of methods \eqref{family of VEM}.
It is worth to stress that the coercivity and continuity constants may depend on $\p$, owing to \eqref{stability}.

Let us set the $H^1$ broken Sobolev seminorm associated with decomposition $\taun$ as:
\begin{equation} \label{H1 broken Sobolev}
\vert \cdot \vert^2_{1,\taun; \Omega} = \sum_{\E \in \taun} \vert \cdot \vert^2_{1,\E}.
\end{equation}
and let us set $S^{\p,-1}(\Omega, \taun)$ the space of piecewise discontinuous polynomials of degree $\p$ over decomposition $\taun$.

Having properties (\textbf{P1}) and (\textbf{P2}), one can prove as in \cite{hpVEMbasic, preprint_hpVEMcorner} an abstract error analysis result.
Given $\u$ and $\un$ the solutions of \eqref{Poisson problem} and \eqref{family of VEM}, respectively, the following holds true:
\begin{equation} \label{abstract error analysis}
\vert \u - \unn \vert_{1,\Omega} \lesssim \frac{\alpha^*(\p)}{\alpha_*(\p)} \left\{ \vert \u - \upi \vert_{1,\taun; \Omega} + \vert \u - \uI \vert_{1,\Omega} + \Fn   \right\} \quad \forall \upi \in S^{\p,-1}(\Omega, \taun),\, \forall \uI \in \Vn,
\end{equation}
where $\Fn$ is the smallest positive constant such that:
\[
\vert (\f,\vnn)_{0,\Omega} - \langle \fn, \vnn \rangle \vert \le \Fn \vert \vnn \vert_{1,\Omega}\quad \forall \vnn \in \Vn.
\]
Thus, being able to study the convergence of the method is equivalent to being able to prove best approximation results by piecewise discontinuous polynomials in $S^{\p,-1}(\Omega, \taun)$
and by functions in $\Vn$, the global VE space defined in \eqref{global VE space}.

In \cite{hpVEMbasic, preprint_hpVEMcorner}, it was proven that for any $\u$, solution of \eqref{Poisson problem}, in $H^{k+1}(\Omega)$, the following $\h\p$ a priori estimate holds true:
\begin{equation} \label{algebraic convergence}
\vert \u - \un \vert_{1,\Omega} \lesssim \frac{\alpha^*(\p)}{\alpha_*(\p)} \frac{\h^{\min(\p,k)}}{\p^k} \vert \u \vert_{k+1,\Omega}
\end{equation}
and that for $\u$ analytic on a proper extension of $\Omega$, the following pure $\p$ approximation result holds true:
\begin{equation} \label{exponential convergence}
\vert \u - \un \vert_{1,\Omega} \lesssim \frac{\alpha^*(\p)}{\alpha_*(\p)} \exp(-b\,\p),
\end{equation}
where $b$ is a positive constant independent of the discretization parameters.

The dependence on $\p$ of the pollution factor $\frac{\alpha^*(\p)}{\alpha_*(\p)}$ is clearly an issue of extreme importance in the $\p$ analysis of the method. The following crude upper bound:
\[
\frac{\alpha^*(\p)}{\alpha_*(\p)} \lesssim \p^{5},
\]
employing a particular stabilization, was proved in \cite[Theorem 4.1]{preprint_hpVEMcorner}.
However, numerical experiments on different stabilizations highlight that the dependence on $\p$ is much milder;
typically, even on nonconvex and nonstar-shaped polygons, one has $\frac{\alpha^*(\p)}{\alpha_*(\p)} \lesssim \p^\beta$, for some $\beta \in (0,1)$;
see \cite{preprint_hpVEMcorner, hpVEMbasic}.

It is worth to stress that such algebraic dependence on the degree of accuracy does not affect the exponential convergence \eqref{exponential convergence} when approximating analytic solution.
Importantly, it is possible to get rid of the dependence on the pollution factor also when approximating solutions presenting corner singularities, employing the $\h\p$ version of the method; see \cite{preprint_hpVEMcorner}.

\subsection{Choices for the stabilization} \label{subsection choices for the stabilization}
Here, we address the issue of choosing a proper stabilization $\SE$ \eqref{stability S} for method \eqref{family of VEM}.
In particular, we define four possible stabilizations which can be found in literature and we discuss their properties.
We recall that $\NdofE$ is the number of dofs of local space $\VnE$.

The first stabilization that we present is:
\begin{equation} \label{classical stabilization}
\SE_1(\unn,\vnn) = \sum_{i=1}^{\NdofE} \dof_i(\unn) \, \dof_i (\vnn),
\end{equation}
which can be regarded as the ``standard'' stabilization of VEM. It was firstly introduced in \cite{VEMvolley} and it can be proved that, for a fixed degree of accuracy $\p$, $\SE_1$ scales like the $H^1$ seminorm.
We highlight that this holds true only in two dimension; for a three dimensional problem one should put a proper scaling factor in front of $\SE_1$, see \cite{preprint_VEM3Dbasic, equivalentprojectorsforVEM}.
The advantage of picking $\SE_1$ as a stabilization is that its implementation is extremely easy.

If one takes into account also variable $\p$, then one has to prove the dependence on $\p$ of the stability constant in \eqref{stability S}.
This was done in \cite{preprint_hpVEMcorner}, where the authors introduced the following \emph{computable} stabilization:
\begin{equation} \label{hpVEMcorner stabilization}
\SE_2(\unn,\vnn) = \frac{\p}{\hE} (\unn, \vnn)_{0,\partial \E} + \frac{\p^2}{\hE^2} (\Pizpmd \unn, \Pizpmd \vnn)_{0,\E}.
\end{equation}
At the current theoretical level, the stabilization constants $c_*(\p)$ and $c^*(\p)$ which appear in \eqref{stability S}
have a dreadful dependence on $\p$; nonetheless, at the practical level, it was observed in \cite{preprint_hpVEMcorner} an extremely mild dependence on $\p$
of $c_*(\p)$ and $c^*(\p)$.

Next, we recall another stabilization which was introduced in \cite{preprint_VEM3Dbasic}. If we denote by $\K_C^\E$ the consistency part of the local stiffness matrix, i.e. the matrix counterpart of the first term in the right-hand side of \eqref{local discrete bilinear form},
then we can define a stabilization $\SE_3$ through its matrix representation $\SEmatthree$ associated with the second term in the right-hand side of \eqref{local discrete bilinear form} as follows.
$\SEmatthree$ is a diagonal matrix, whose $i$-th entry is given by:
\begin{equation} \label{stabilization Dassi}
(\SEmatthree)_{i,i} = \max \left(1, (\mathbf{K}_C^\E)_{i,i} \right).
\end{equation}
For the sake of clarity, the complete local stiffness matrix $\mathbf{K}_n^\E$ reads:
\[
\K_\p^\E = \K_C^\E  + \K_{S_3}^\E = \K _C^\E + (\Id - \Pinablamat)^T \cdot \SEmatthree \cdot (\Id - \Pinablamat),
\]
where $\Pinablamat$ is the matrix representing the action of operator $\Pinablap$ and $\SEmatthree$ is defined in \eqref{stabilization Dassi}.

In \cite{preprint_VEM3Dbasic}, numerical experiments, along with a heuristic motivation, show that  choice \eqref{stabilization Dassi} entails a better behaviour of the method for high $\p$, when approximating a 3D Poisson problem.

Finally, we also consider a fourth stabilization. Let $\Ndofbndr$ be the number of boundary degrees of freedom of $\VnE$.
Then we set:
\begin{equation} \label{stabilization only boundary}
\SE_4(\unn, \vnn) = \sum_{i=1}^{\Ndofbndr} \dof_i(\unn) \dof_i(\vnn),
\end{equation}
which is basically stabilization $\SE_1$ without the contribution of internal dofs.

We point out that other choices of stabilization can be found in literature, but we prefer to investigate the four presented here in order to avoid an overwhelming set of numerical experiments.

We also emphasize that the first three stabilizations above hinge upon the choice of the polynomial basis dual to internal moments \eqref{internal moments}.

\subsection{Choices for the basis} \label{subsection choices for the basis}
In this section, we discuss three possible choices of the local VE basis. More precisely, we consider internal moments \eqref{internal moments} taken with respect to three different polynomial bases
$\{ \q_\alpha^i \} _{\alpha=1}^{\dim(\mathbb P_{\p-2}(\E))}$, $i=1,2,3$, of $\mathbb P_{\p-2}(\E)$, thus modifying the definition of the VE bubble functions.

The hope is that a proper choice of the polynomial basis dual to internal moments \eqref{internal moments}, and therefore of internal VEM basis elements, entails a better conditioning for the stiffness matrix.

Henceforth, we will employ, with an abuse of notation, the natural bijection $\mathbb N^2 \leftrightarrow  \mathbb N_0 $ given by:
\begin{equation} \label{natural bijection}
(0,0) \leftrightarrow 1,\;\;\;\; (1,0) \leftrightarrow {2},\;\;\;\; (0,1) \leftrightarrow 3,\;\;\;\; (2,0) \leftrightarrow  4,\;\;\;\; (1,1) \leftrightarrow  5,\;\;\;\; (0,2) \leftrightarrow  6,\;\;\;\;\dots
\end{equation}
We will also write occasionally:
\[
\qalphabasisi \azpmd \quad \text{instead of $\{ \q_\alpha^i \} _{\alpha=1}^{\dim(\mathbb P_{\p-2}(\E))}$}, \quad \text{for $i=1,2,3$}.
\]

The first choice of the polynomial basis is the ``standard'' one, i.e. the one which is used in the majority of VEM literature, since, at the implementation level, is the most convenient.
Given $\xEbold = (\xE, \yE)$ the barycenter of $\E$, we set:
\begin{equation} \label{first choice basis}
\qalphau (\xbold) = \left( \frac{\xbold - \xEbold}{\hE}  \right)^{\boldalpha} = \left( \frac{x-\xE}{\hE} \right)^{\alpha_1}  \left( \frac{y-\yE}{\hE} \right)^{\alpha_2} \quad \forall \boldalpha = (\alpha_1,\alpha_2) \in \mathbb N^2,\quad \vert \boldalpha \vert = 0,\dots, \p-2.
\end{equation}
Although choice \eqref{first choice basis} is very suitable from the computational point of view, it has bad effects on the condition number of the stiffness matrix for high local degrees of accuracy $\p$,
see \cite{hpVEMbasic}, and in presence of bad-shaped polygons, see \cite{BerroneBorio_orthogonalVEM} and the references therein.

For this reason, we suggest two possible modifications which rely on orthogonalizing processes of $\qalphabasisu \azpmd$ with respect to the $L^2$ norm on polygon $\E$.

The first modification, which allows to construct an $L^2(\E)$ orthonormal basis $\qalphabasisd \azpmd$, is based on the stable Gram-Schmidt orthonormalization process presented in \cite{BassiBottiColomboDipietroTesini}.
We point out that basis $\qalphabasisd \azpmd$ was firstly introduced in the context of $\p$-VEM multigrid algorithm for the construction of the multrigrid scheme, see \cite{pVEMmultigrid}.

While the implementation of VEM with basis $\qalphabasisu \azpmd$ is well-known (see \cite{hitchhikersguideVEM}), no details were given in \cite{pVEMmultigrid} for the implementation of the method using basis $\qalphabasisd \azpmd$.
Therefore, we address this issue in Appendix \ref{section hitchhikers}.

The third choice of the polynomial basis is slightly inspired by an orthonormalizing procedure used in \cite{BerroneBorio_orthogonalVEM}.
In order to present the third basis $\qalphabasist \azpmd$, we set matrix:
\[
\H _{\boldalpha, \boldbeta} = (\qalphau, \q_{\boldbeta}^1)_{0,\E} \quad \forall \boldalpha,\, \boldbeta \in \mathbb N^2, \; \vert \boldalpha \vert,\, \vert \boldbeta \vert =0,\dots, \p-2.
\]
Given $n_\ell = \dim (\mathbb P_\ell(\E))$, we fix $\qalphat = \qalphau$ if $\alpha= 1$ and we decompose matrix $\H$ into blocks as:
\[
\H = 
\bordermatrix{ 			    & \scriptstyle {1} & \scriptstyle{\npmd-1} \cr
                  \scriptstyle 1               & \H_{1,1} & \H_{1,2} \cr
                  \scriptstyle{\npmd-1} &\H_{2,1} & \H_{2,2} \cr},
\]
we diagonalize matrix $\H_{2,2}$ and we get:
\[
\mathbf V^t \cdot \H_{2,2} \cdot \mathbf V = \D \Longrightarrow  (\mathbf V \cdot \D^{-\frac{1}{2}})^t \cdot \H_{2,2} \cdot (\mathbf V \cdot \D^{-\frac{1}{2}}) = \Id.
\]
This entails that matrix $\mathbf V \cdot \D^{-\frac{1}{2}}$ contains the coefficients which orthonormalize $\qalphabasisu \aupmd$, the monomial basis of $\mathbb P_{\p-2}(\E) / \mathbb R$. Therefore, one has:
\begin{equation} \label{third choice basis}
\qalphat (\xbold) = 
\begin{cases}
\qalphau & \text{if } \vert \boldalpha \vert=0\\
\sum_{\vert \boldbeta \vert=1}^{\p-2} (\mathbf V \cdot \D^{-\frac{1}{2}}) ^t _{\boldalpha, \boldbeta} \;\qbeta^1(\xbold) & \text{if } \vert \boldalpha \vert =1,\dots, \p-2\\
\end{cases}.
\end{equation}
It is worth to stress that here the orthonormalizing process is performed with a different target with respect to what done in \cite{BerroneBorio_orthogonalVEM};
in fact, there, the method is built employing the canonical bases computed taking moments with respect to scaled monomials $\qalphabasisu \azpmd$;
however, the projectors onto polynomial spaces, i.e. $\Pizpmd$ and $\Pinablap$ defined in \eqref{L2 orthogonality} and \eqref{H1 orthogonality}, respectively,
are computed expanding with respect to polynomial basis $\qalphabasist \azpmd$ defined in \eqref{third choice basis}.
Here, we define in addition internal moments with respect to basis $\qalphabasist \azpmd$.

The implementation issues regarding this new basis are not explicated in this paper, but are similar to those of basis $\qalphabasisd \azpmd$, which are shown in Appendix \ref{section hitchhikers}.

In summary, we presented three choices for the polynomial basis dual to internal moments \eqref{internal moments}. One is of easy implementation, but it may be the cause of a high condition number for the stiffness matrix
when using high-order methods or in presence of ``bad-shaped'' polygons.
The other two bases are obtained by two distinct orthonormalization processes; their performances with respect to the condition number are investigated in Section \ref{section numerical results}.
What we can anticipate is that they outclass the performances of their counterpart using the standard basis $\qalphabasisu \azpmd$ in the two situations above mentioned.

A \emph{heuristic} reason for this fact is the following. If, for each element $\E$, the local Virtual Element Space $\VnE$ were a space consisting of polynomials only,
then picking internal moments \eqref{internal moments} with respect to an $L^2(\E)$ orthonormal polynomial basis would automatically entail
that the local canonical basis is made of polynomials and contain a subset of $L^2(\E)$ orthonormal polynomials spanning $\mathbb P_{\p-2}(\E)$;
it is well-known in the theory of Spectral Elements, see e.g. \cite{pasquettirapetti2004spectral, canuto2012spectral}, that
employing $L^2$ orthonormal canonical basis damp the condition number of the stiffness matrix when increasing the polynomial degree.

Nonetheless, local Virtual Element Spaces does not contain polynomials only, but also other functions needed for prescribing $H^1$ conformity.
As stated in \eqref{dimension of subspaces}, the dimension of the subspace of nonpolynomial functions is, in terms of the degree of accuracy, asymptotically smaller than the dimension of the subspace of polynomial functions.
Therefore, \emph{in a very rough sense}, employing $L^2$ orthonormal polynomials in the definition of internal moments entails a sort of \emph{partial ``orthonormalization''} of the local canonical basis.


Before concluding this section, we associate to bases $\qalphabasisi \azpmd$ the sets of dofs $\{ \dof^i_j\}_{j=1}^{\NdofE}$ and the canonical bases $\{\varphi_j^i\}_{j=1}^{\NdofE}$,  for all $i=1,2,3$.
This notation will be instrumental in Appendix \ref{section hitchhikers}.

\subsection{Stabilizations and bases: the effects on the method} \label{subsection stab and bases effects on the method}
Having presented in the two foregoing sections various choices of stabilizations and canonical bases,
we want here to highlight the effects of such choices on the method and on the ill-conditioning of the stiffness matrix.

\paragraph*{Stabilization.}
The choice of the stabilization has two effects. The first one is related to the convergence of the method
since nonproperly tailored choices of the stabilization automatically entail higher pollution factor $\frac{\alpha^*(\p)}{\alpha_*(\p)}$, see \eqref{exponential convergence} and especially \eqref{algebraic convergence}.
Secondly, since the stabilization appears in the discrete bilinear form of the method, see \eqref{local discrete bilinear form}, there is also an effect on the condition number of the stiffness matrix.

\paragraph*{Canonical basis.}
The choice of the canonical basis has also two effects.
Firstly, it has an impact on the condition number of the global stiffness matrix, simply because by changing the basis automatically the entries of the stiffness matrix modify.
Secondly, by picking different canonical bases, one also changes the definition of the stabilization;
as an example, if we fix stabilization $\SE_1$ defined in \eqref{classical stabilization} and we apply it to functions in the Virtual Element space, then we get in general different values,
since the definition of the internal degrees of freedom vary depending on the choice of the basis; in particular, one also modify the behaviour of the pollution factor $\frac{\alpha^*(\p)}{\alpha_*(\p)}$.

\section {Numerical results} \label{section numerical results}
In this section we present some numerical experiments in which we compare the performances of the stabilizations and polynomial bases introduced in Sections
\ref{subsection choices for the stabilization} and \ref{subsection choices for the basis}, respectively.

More precisely, we investigate the behaviour and the related effects of the condition number in two critical situations.


In Section \ref{subsection numerical results: increasing p}, we investigate the behaviour of the condition number of the $\p$ version of VEM and the effects on the linear solver used for the resolution of the associated linear system.
We are also interested in the behaviour of the condition number when varying the stabilization and the polynomial basis dual to internal moments \eqref{internal moments} in presence of a sequence of ``bad shaped'' polygons
(collapsing bulk, collapsing edges\dots);
this is probed in Sections \ref{subsection numerical results: collapsing polygon} and \ref{subsection numerical results: hanging nodes}.
The condition number is computed as the ratio between the maximum and the minimum (nonzero) eigenvalue of the stiffness matrix.

\subsection{Numerical results: the $\p$ version of VEM} \label{subsection numerical results: increasing p}
\begin{figure}  [h]
\centering
\subfigure {\includegraphics [angle=0, width=0.32\textwidth]{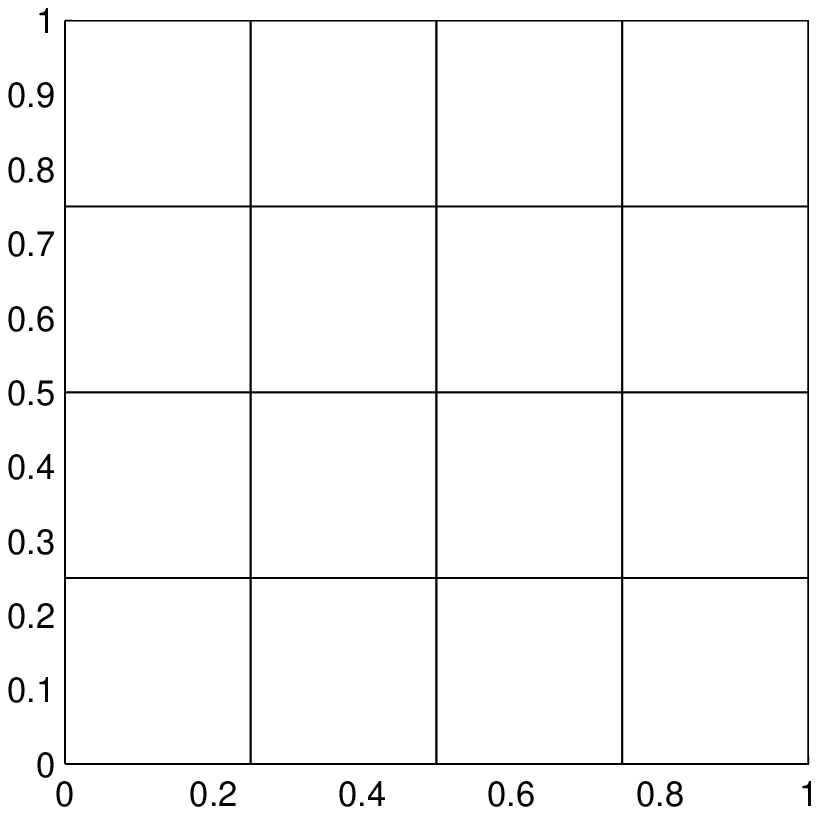}}
\subfigure {\includegraphics [angle=0, width=0.32\textwidth]{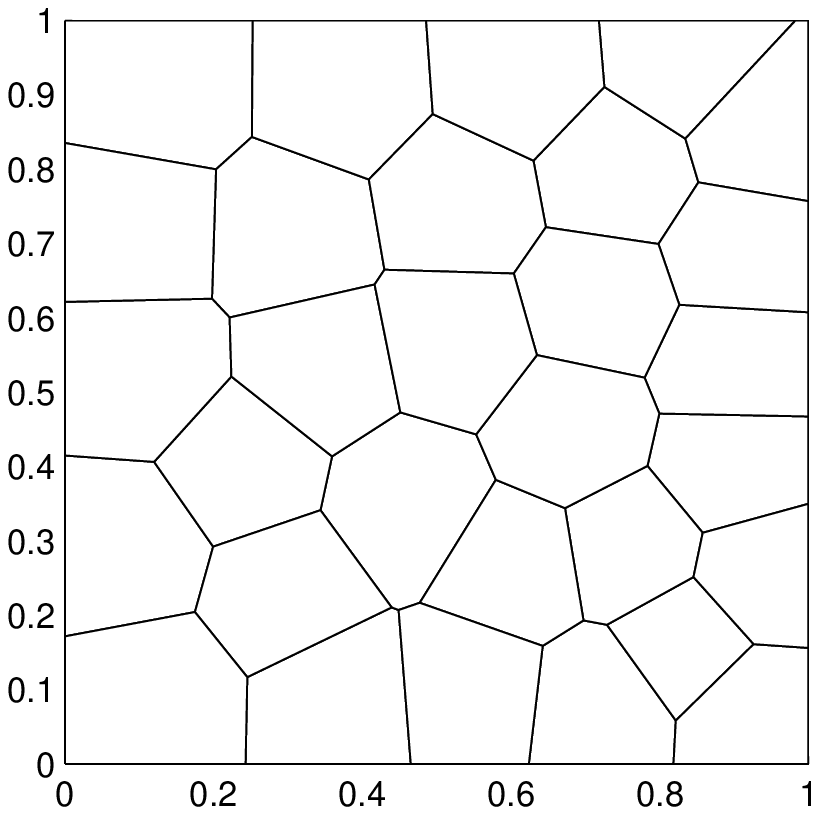}}
\subfigure {\includegraphics [angle=0, width=0.32\textwidth]{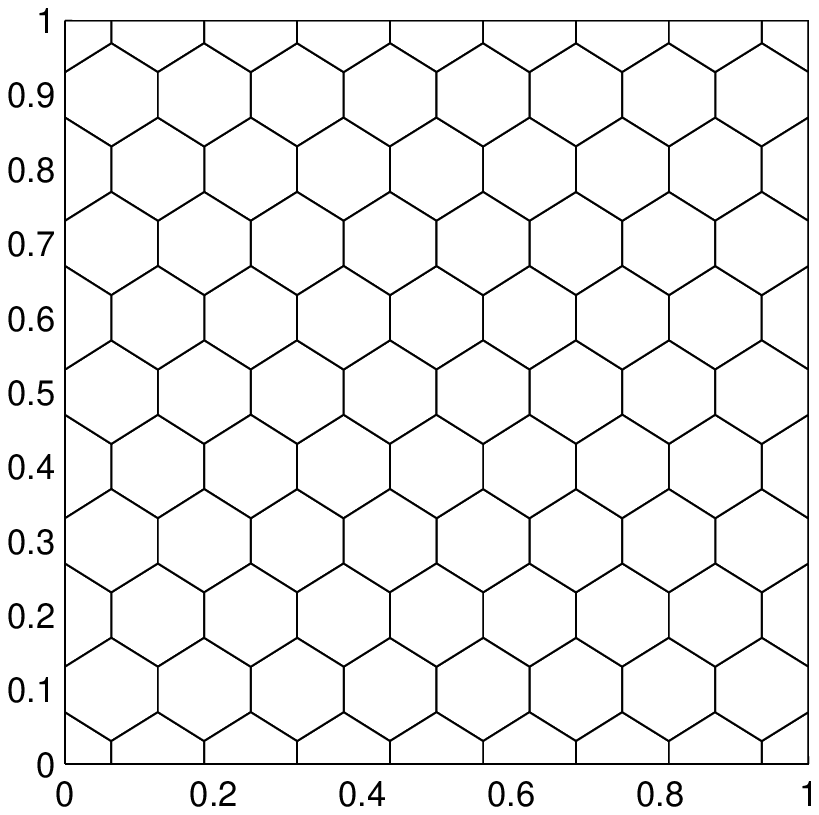}}
\caption{Left: square mesh. Center: Voronoi-Lloyd mesh. Right: regular-hexagonal mesh.} \label{figure three meshes}
\end{figure}

Let us consider the three meshes depicted in Figure \ref{figure three meshes}.
We investigate in this section the behaviour of the condition number of the stiffness matrix associated with method \eqref{family of VEM} by keeping fixed the meshes of Figure \ref{figure three meshes} and by increasing $\p$.
For the purpose, we modify the choice of the stabilization and the choice of the polynomial basis dual to internal moments \eqref{internal moments}.

\begin{figure}  [h]
\centering
\subfigure {\includegraphics [angle=0, width=0.32\textwidth]{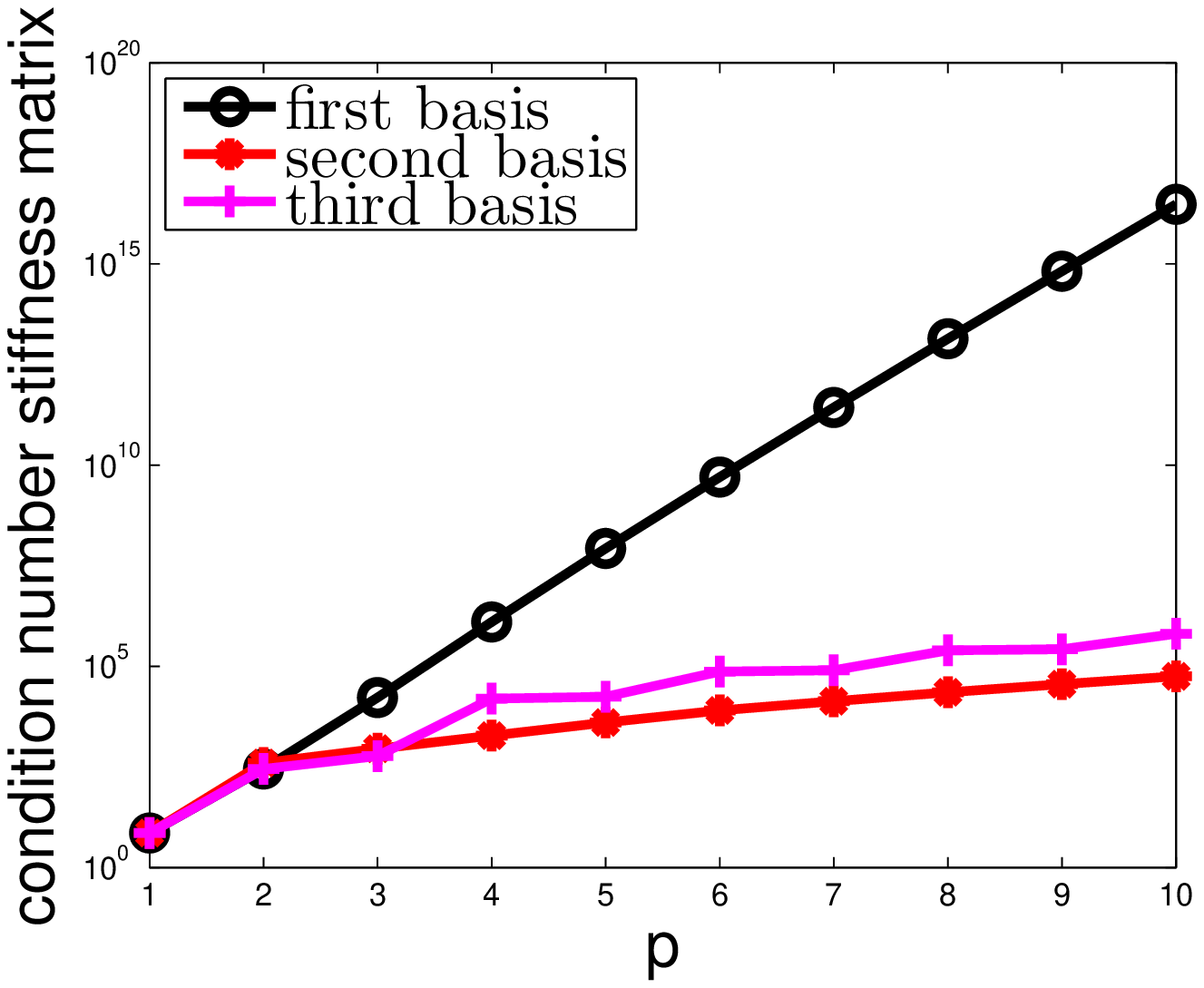}}
\subfigure {\includegraphics [angle=0, width=0.32\textwidth]{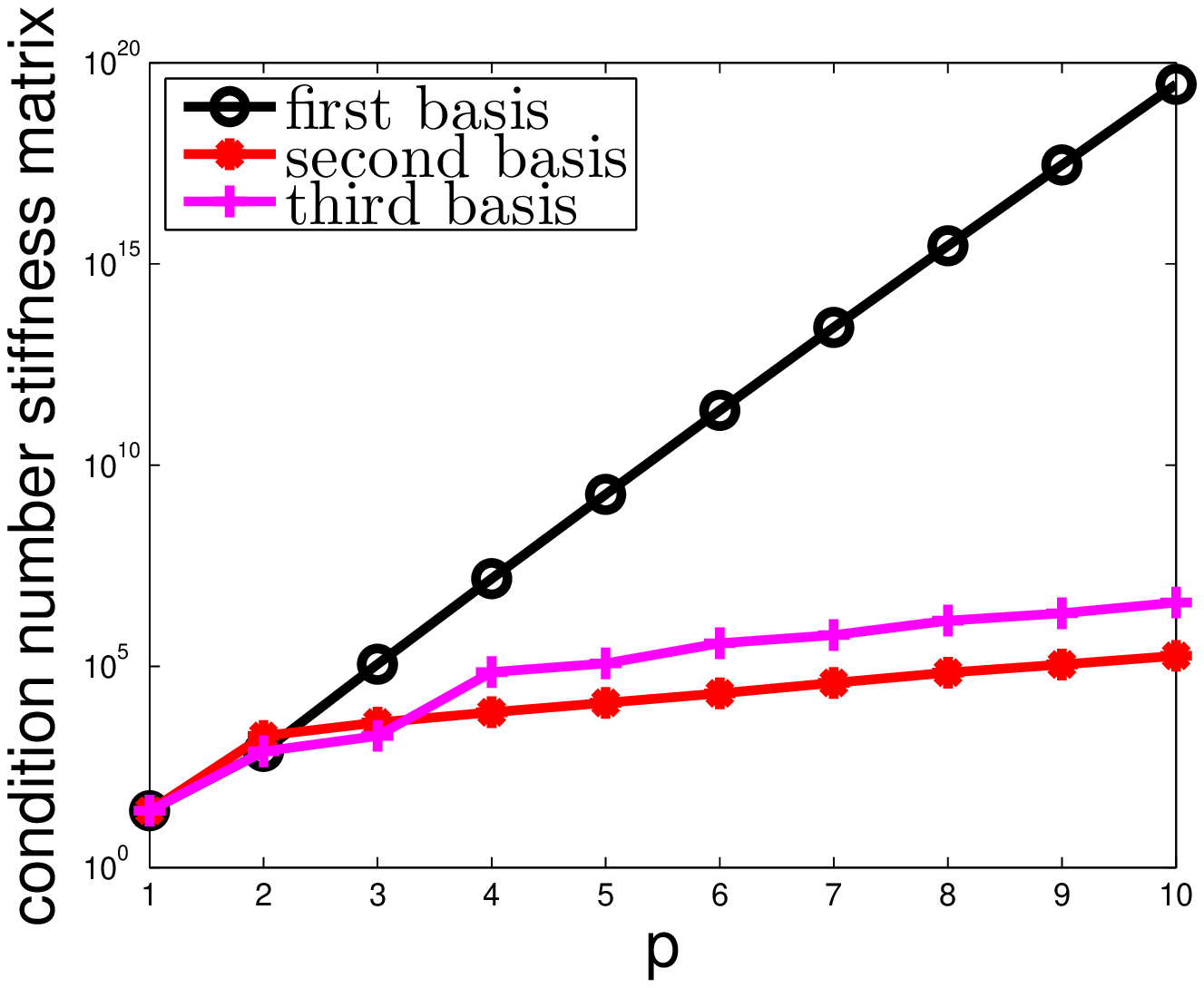}}
\subfigure {\includegraphics [angle=0, width=0.32\textwidth]{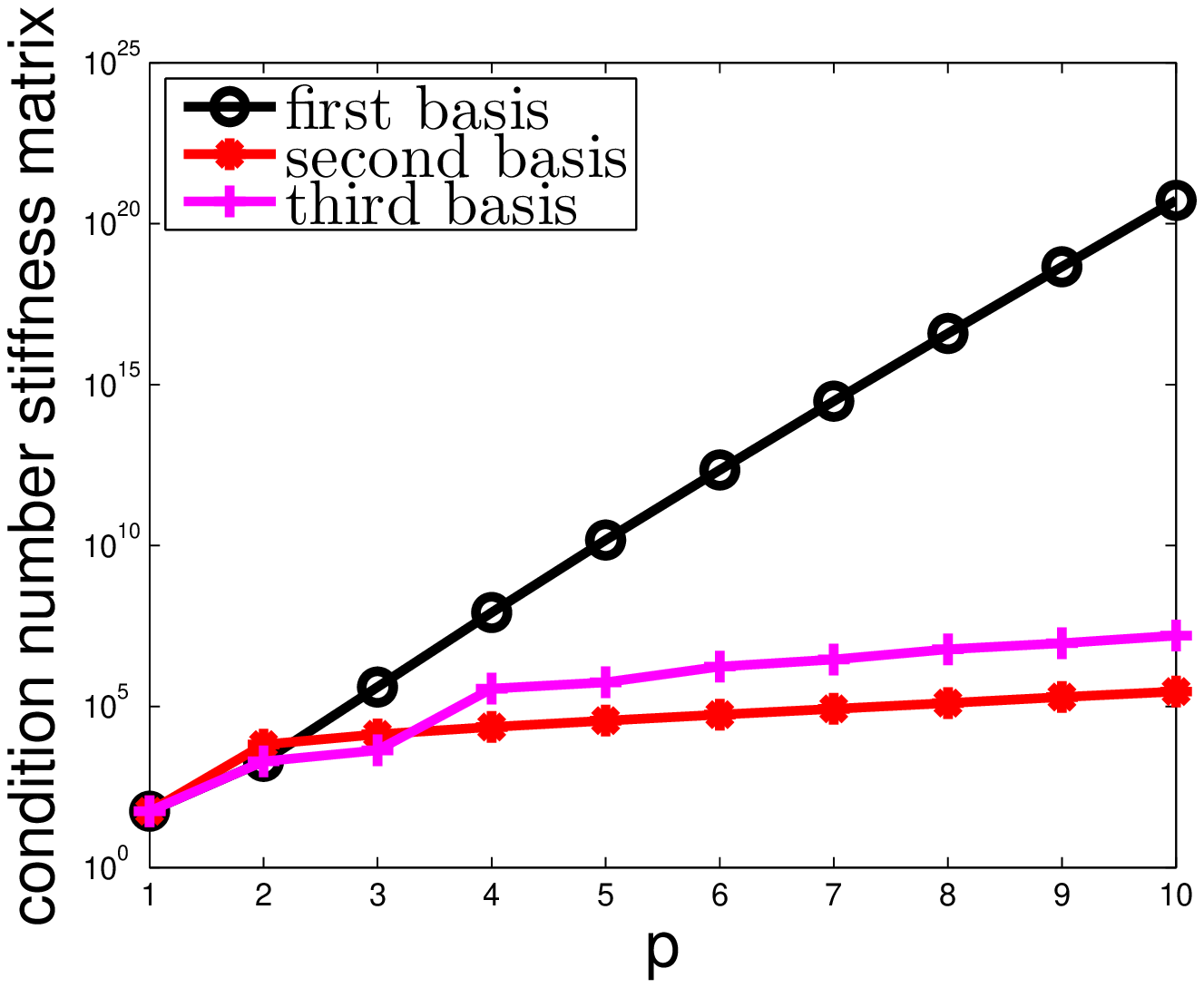}}
\caption{Condition number of the VEM stiffness matrix in terms of $\p$.
The stabilization is fixed and equal to $\SE_1$ \eqref{classical stabilization}. We compare the behaviour in terms of the three polynomial bases presented in Section \ref{subsection choices for the basis}.
Left: square mesh. Center: Voronoi-Lloyd mesh. Right: regular-hexagonal mesh.} \label{figure p cond classical stab analytic solution}
\end{figure}

In Figure \ref{figure p cond classical stab analytic solution}, we depict the behaviour of the condition number by fixing the stabilization to be $\SE_1$ presented in \eqref{classical stabilization} and we consider
the three polynomial bases introduced in Section \ref{subsection choices for the basis}.
For all the three meshes, the basis $\qalphabasisd \azpmd$ shows the best performances, whereas the standard monomial basis $\qalphabasisu \azpmd$ shows the worst results.

From Figure \ref{figure p cond classical stab analytic solution}, it is also clear that basis $\qalphabasisu \azpmd$ entails an \emph{exponential} growth of the condition number in terms of $\p$.

Furthermore, employing $\qalphabasisi \azpmd$, $i=2,3$, suggests instead an algebraic growth of the condition number in terms of $\p$. A polynomial fitting yields:
\begin{equation} \label{algebraic cond}
\text{cond}(\mathbf K_\p) \approx a \, \p^b \quad \quad \quad \text{with }
a = \begin{cases}
130.4 & \text{if } i=2\\
131.7 & \text{if } i=3\\
\end{cases},\quad
b= \begin{cases}
3.344 & \text{if } i=2\\
3.371 & \text{if } i=3\\
\end{cases}.
\end{equation}
This behaviour is extremely interesting since it is well-known, see e.g. \cite{pasquettirapetti2004spectral}, that the growth in terms of $\p$ of the condition number in triangular Spectral Elements with nodal bases is of the following sort:
\begin{equation} \label{cond SE on trinagles}
\text{cond}(\mathbf K_\p) \approx a \, \p^b \quad \quad \quad \text{with } b = 4,\quad \text{for some } a>0.
\end{equation}

We want now to understand how much the ill-conditioning pollutes the convergence of the error:
\begin{equation} \label{error of the method}
\vert \u - \Pinablap \un \vert_{1,\taun} := \sqrt{\sum_{\E \in \taun} \vert \u - \Pinablap \un \vert^2_{1,\E}},
\end{equation}
where we recall that $\u$ is the solution of \eqref{Poisson problem} and $\un$ is its VEM approximation.

\begin{figure}  [h]
\centering
\subfigure {\includegraphics [angle=0, width=0.32\textwidth]{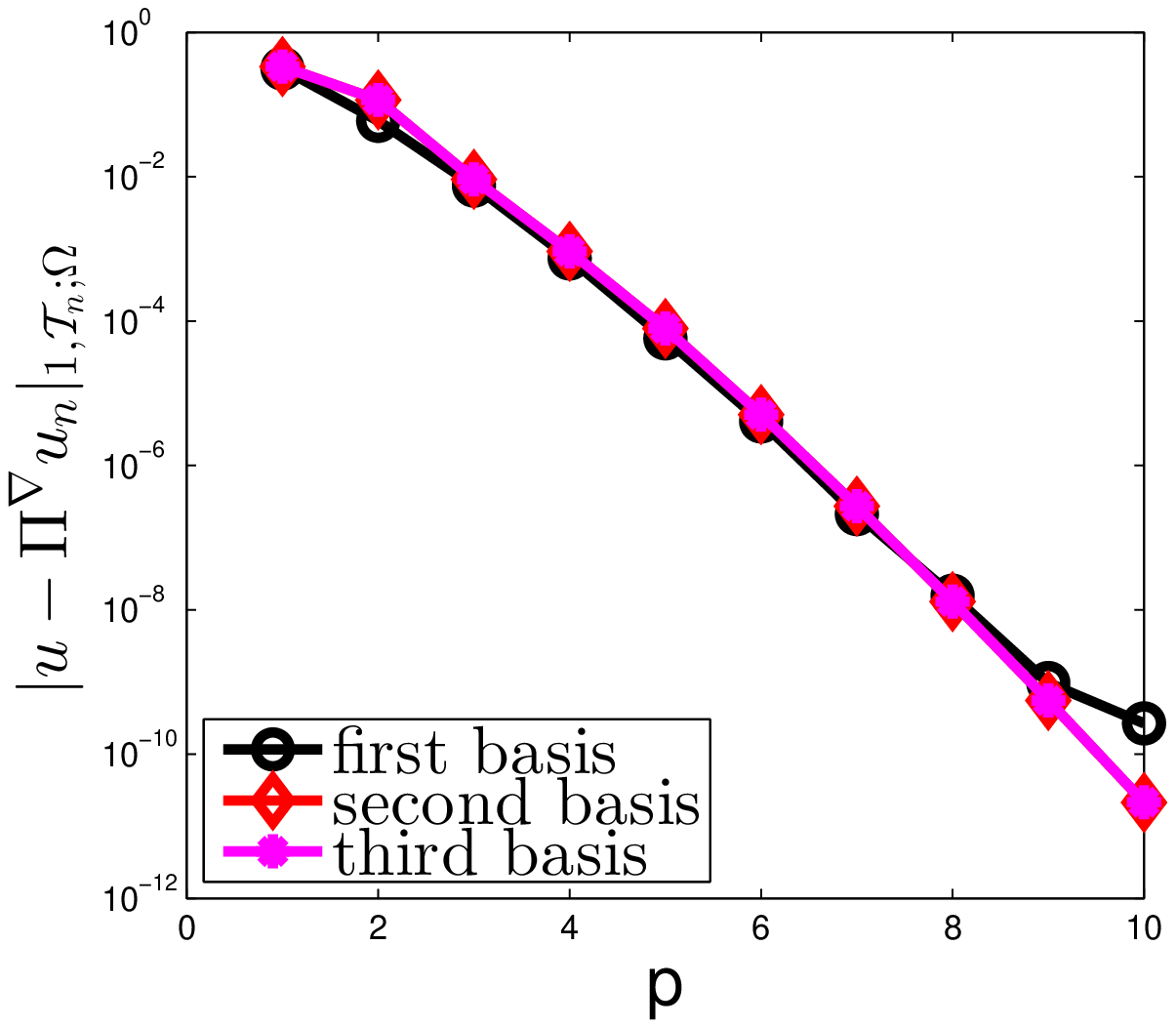}}
\subfigure {\includegraphics [angle=0, width=0.32\textwidth]{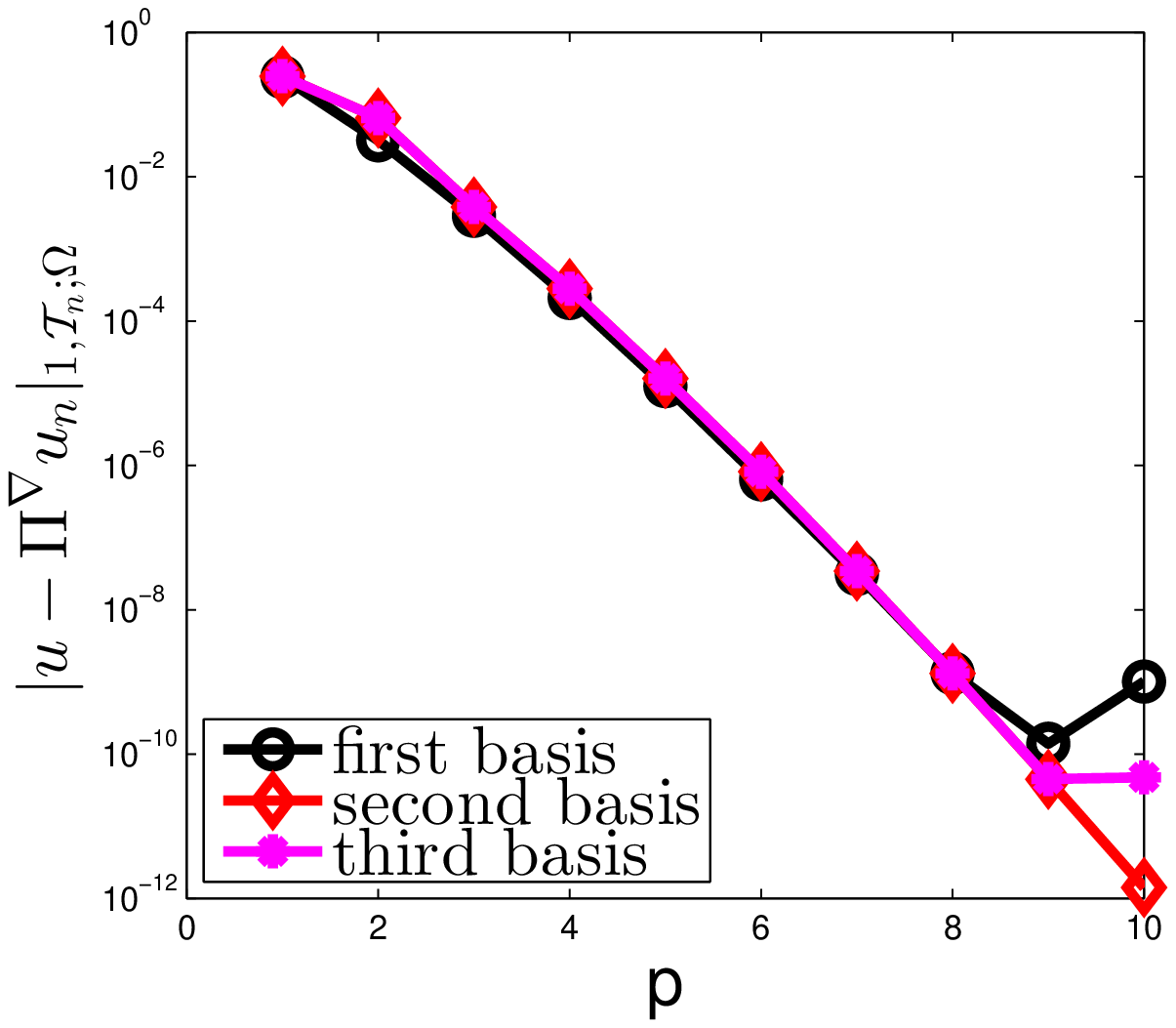}}
\subfigure {\includegraphics [angle=0, width=0.32\textwidth]{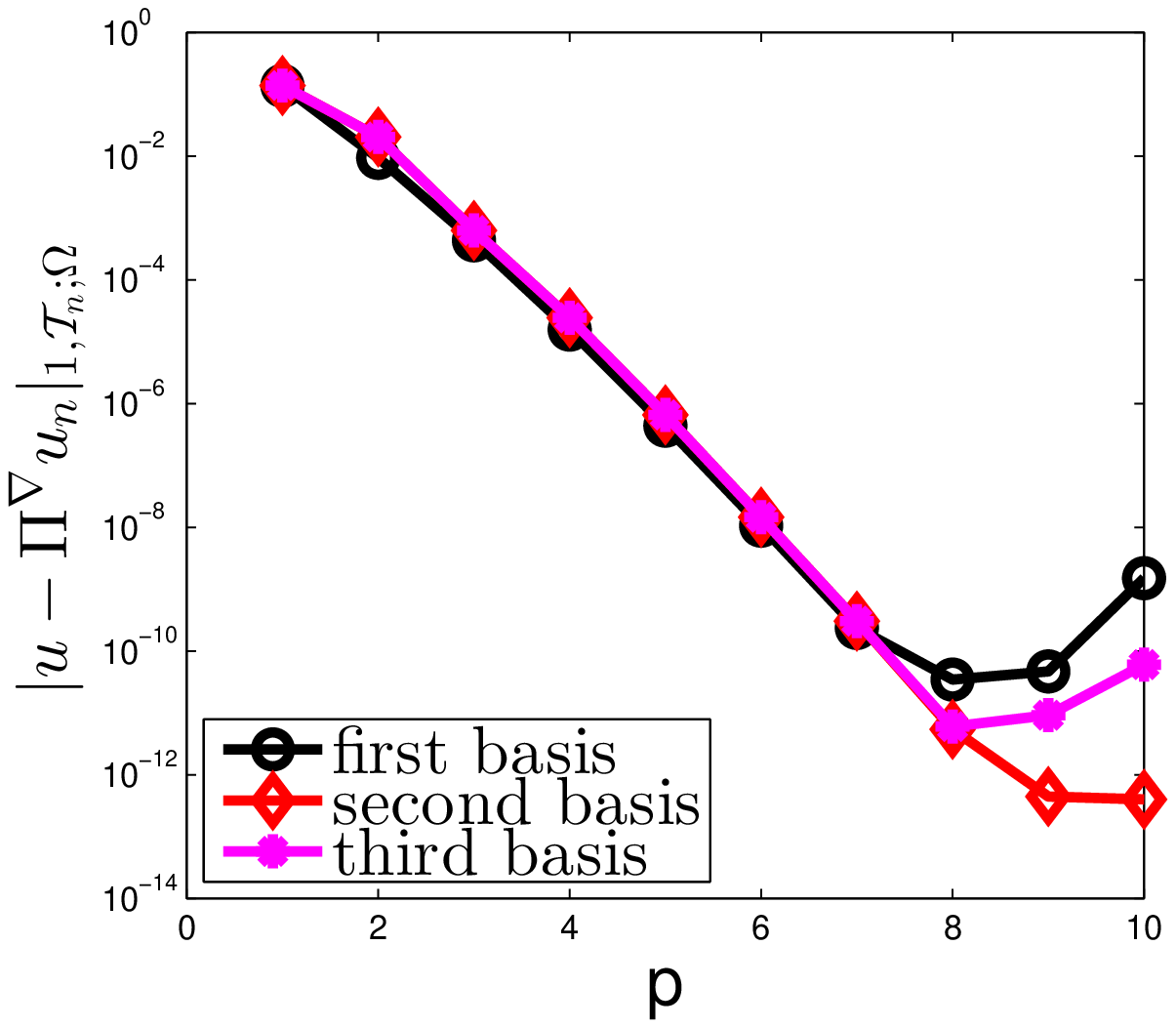}}
\caption{Error $\vert \u - \Pinablap \un \vert_{1,\taun}$ with exact solution given in \eqref{analytic solution}.
The stabilization is fixed and equal to $\SE_1$ \eqref{classical stabilization}. We compare the behaviour in terms of the three polynomial bases presented in Section \ref{subsection choices for the basis}.
Left: square mesh. Center: Voronoi-Lloyd mesh. Right: regular-hexagonal mesh.} \label{figure p error analytic solution various bases}
\end{figure}

For the purpose, we consider a test case with analytic solution:
\begin{equation} \label{analytic solution}
\u(x,y) = \sin(\pi x) \sin(\pi y),
\end{equation}
for which we know that the method converges exponentially, see \eqref{exponential convergence}. In Figure \ref{figure p error analytic solution various bases},
we compare the errors \eqref{error of the method} using the three meshes in Figure \ref{figure three meshes} (always using $\SE_1$ \eqref{classical stabilization} as a stabilization) and comparing the three bases $\qalphabasisi \azpmd$, $i=1,2,3$.
We observe that, due to the ill-conditioning of basis $\qalphabasisu \azpmd$ for high values of $\p$, the linear solver of the system (namely the one associated with the $\backslash$  command of \textsc{MATLAB}) does not work properly.
For this reason, we highly recommend to use basis $\qalphabasisd \azpmd$  in lieu of basis $\qalphabasisu \azpmd$ when approximating with high-order VEM.

In order to understand better ``how much the (linear) solver fails'' when solving the system arising from \eqref{family of VEM}, we consider as an exact solution:
\begin{equation} \label{patch test}
\u(x,y) = 1-x-y,
\end{equation}
which, owing to polynomial consistency assumption \eqref{consistency}, should be approximated exactly by the VEM.

In exact-arithmetic one would expect the error \eqref{error of the method} to vanish, while in floating-point arithmetic the error is not zero but grows along with the condition number of the stiffness matrix.

\begin{figure}  [h]
\centering
\subfigure {\includegraphics [angle=0, width=0.32\textwidth]{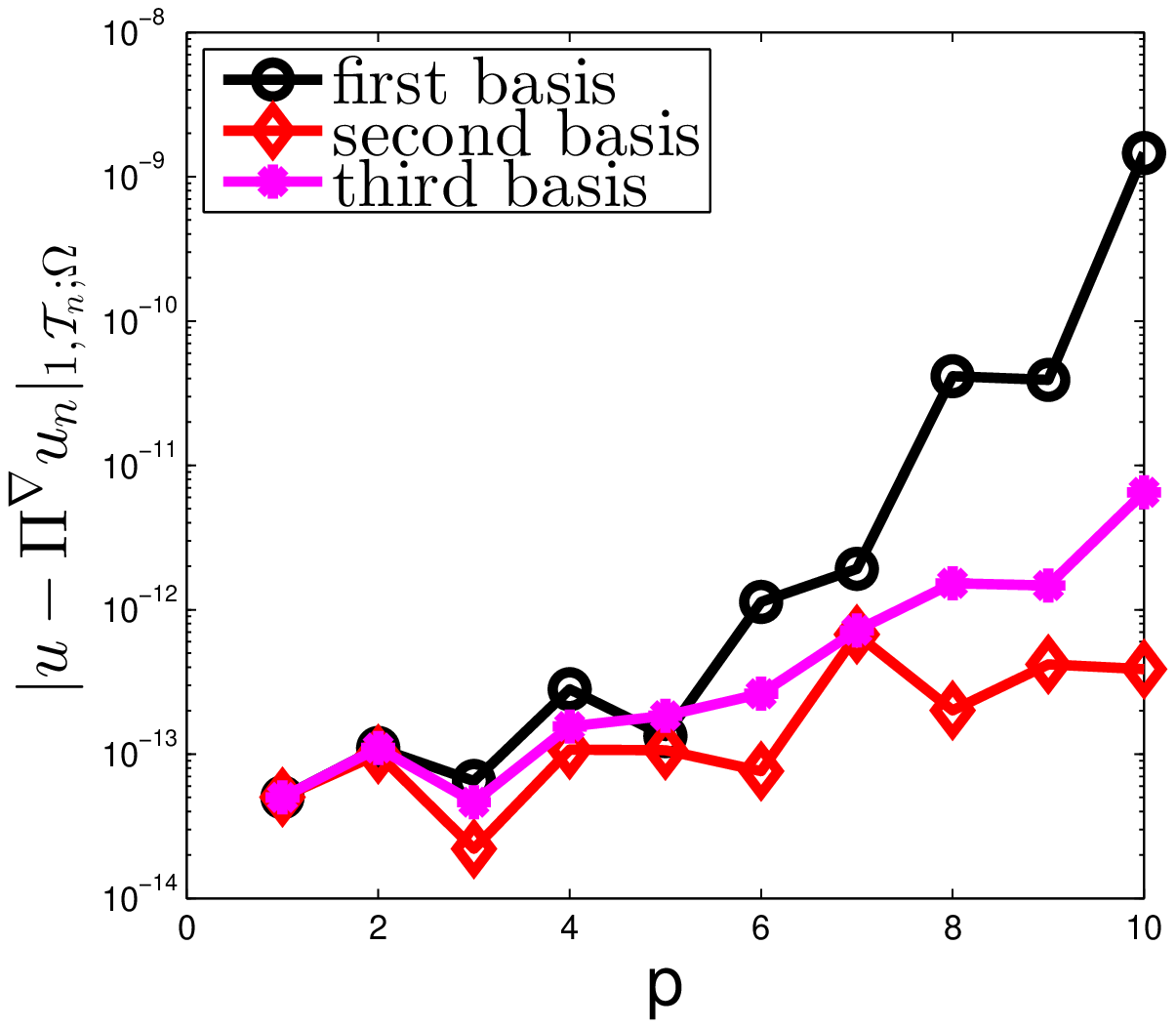}}
\subfigure {\includegraphics [angle=0, width=0.32\textwidth]{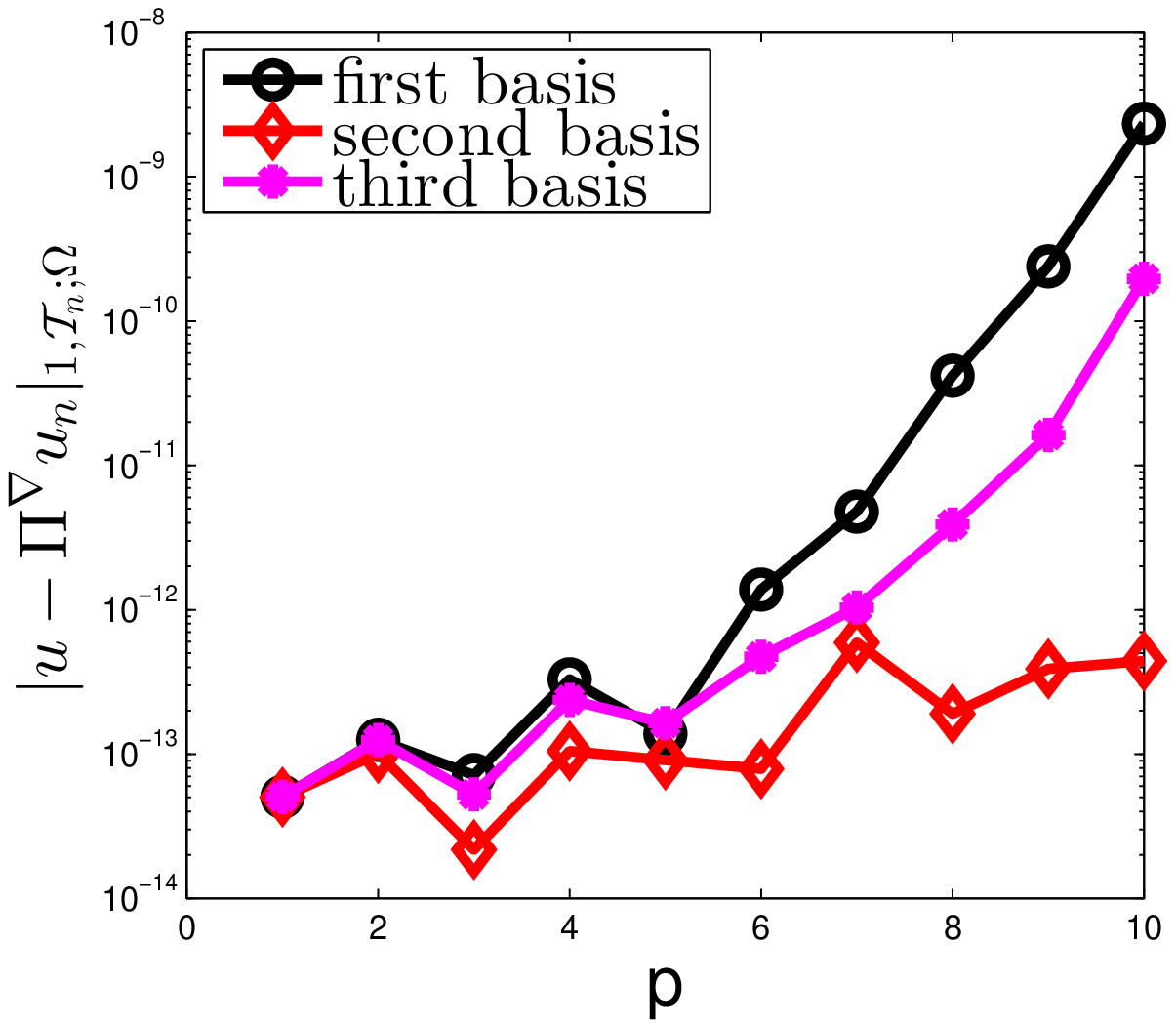}}
\subfigure {\includegraphics [angle=0, width=0.32\textwidth]{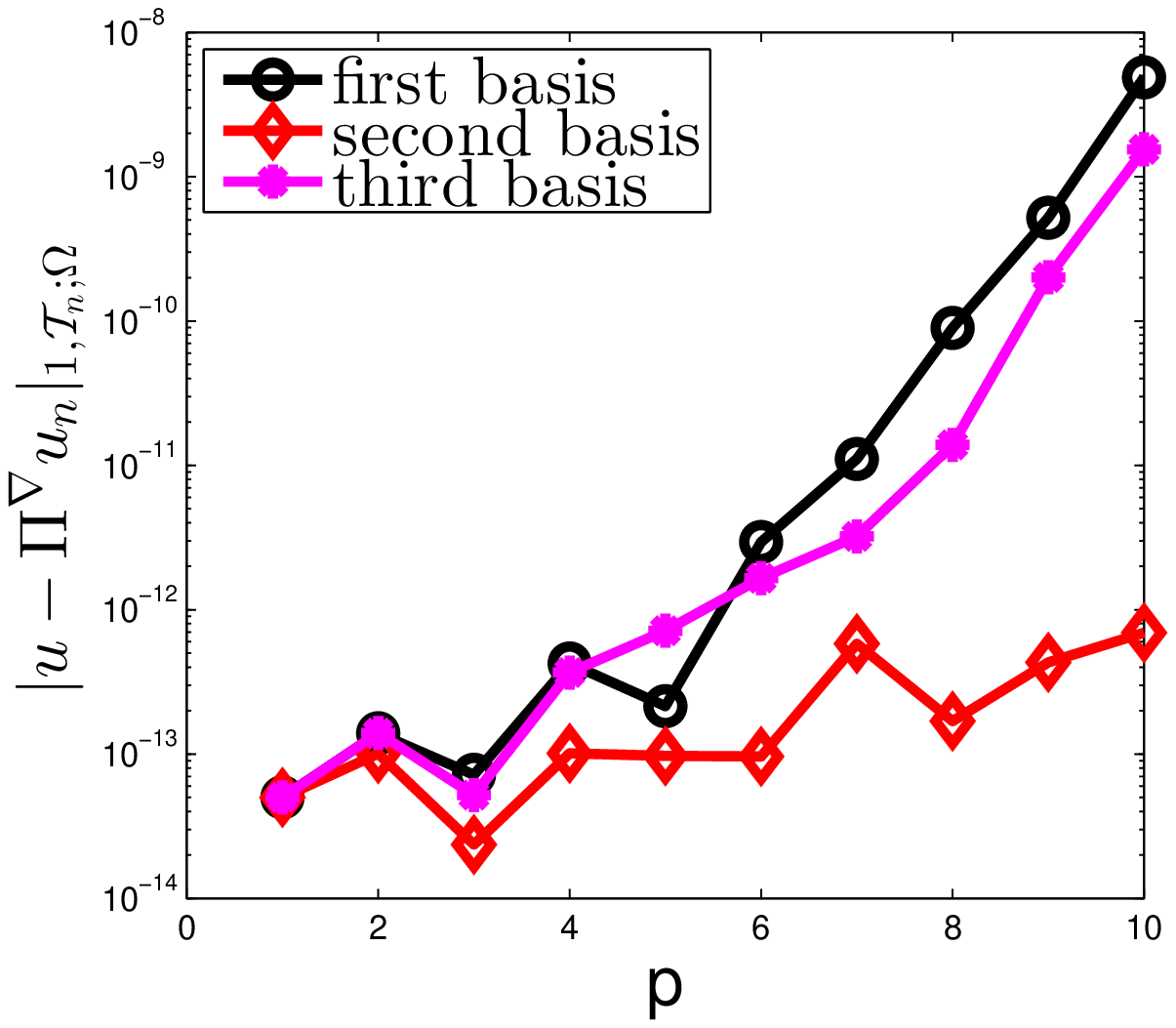}}
\caption{Error $\vert \u - \Pinablap \un \vert_{1,\taun}$ with exact solution given in \eqref{patch test}.
The stabilization is fixed and equal to $\SE_1$ \eqref{classical stabilization}. We compare the behaviour in terms of the three polynomial bases presented in Section \ref{subsection choices for the basis}.
Left: square mesh. Center: Voronoi-Lloyd mesh. Right: regular-hexagonal mesh.} \label{figure p error patch test various basis}
\end{figure}

In Figure \ref{figure p error patch test various basis}, we compare error \eqref{error of the method} using the three meshes in Figure \ref{figure three meshes}, using $\SE_1$ \eqref{classical stabilization} as a stabilization and the
three bases $\qalphabasisi \azpmd$, $i=1,2,3$.
The behaviour of basis $\qalphabasisd \azpmd$ is again superior to the other two bases.
More precisely, employing basis $\qalphabasisu \azpmd$ has a large effect on the error for high degrees of accuracy $\p$.

\begin{figure}  [h]
\centering
\subfigure {\includegraphics [angle=0, width=0.32\textwidth]{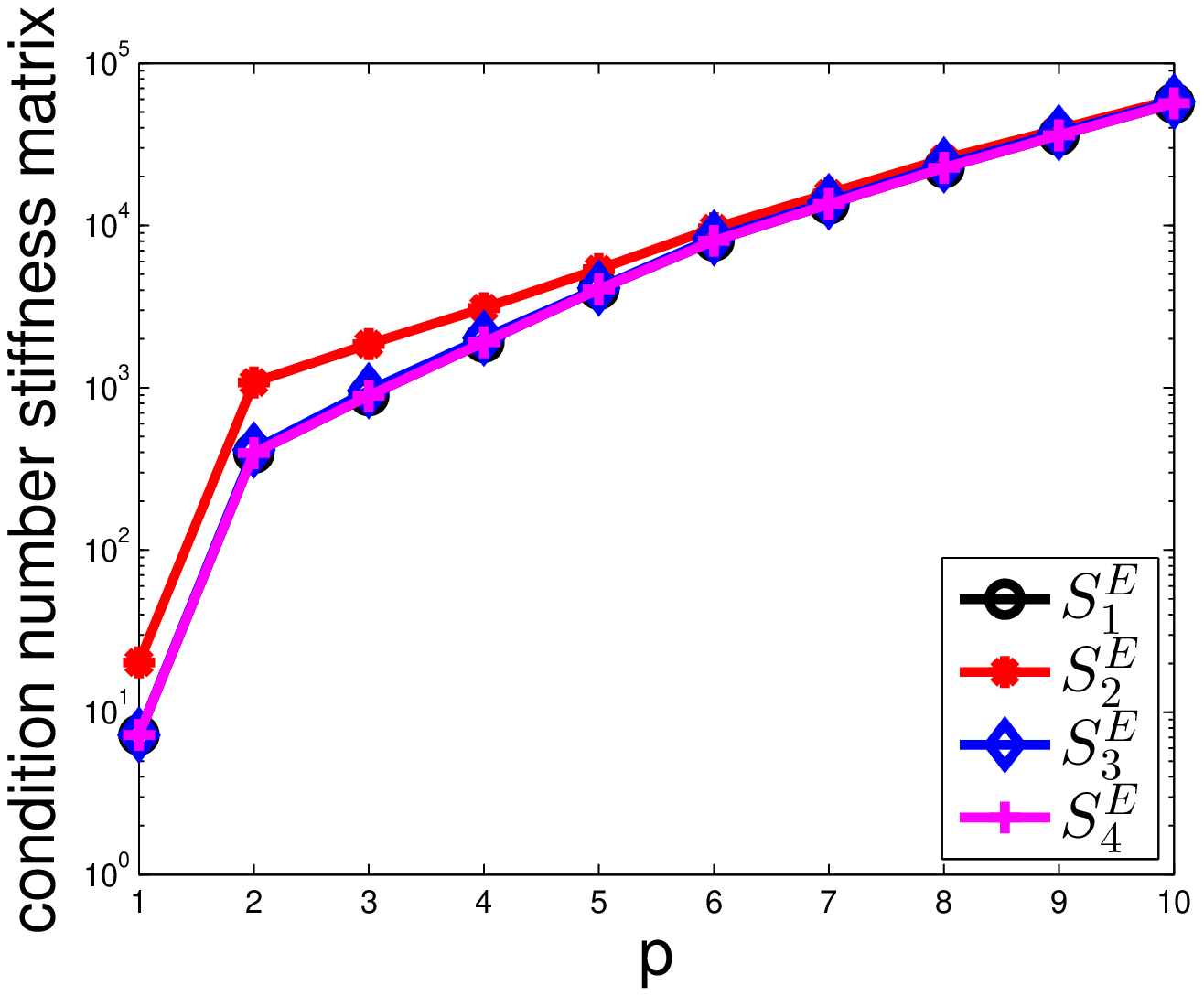}}
\subfigure {\includegraphics [angle=0, width=0.32\textwidth]{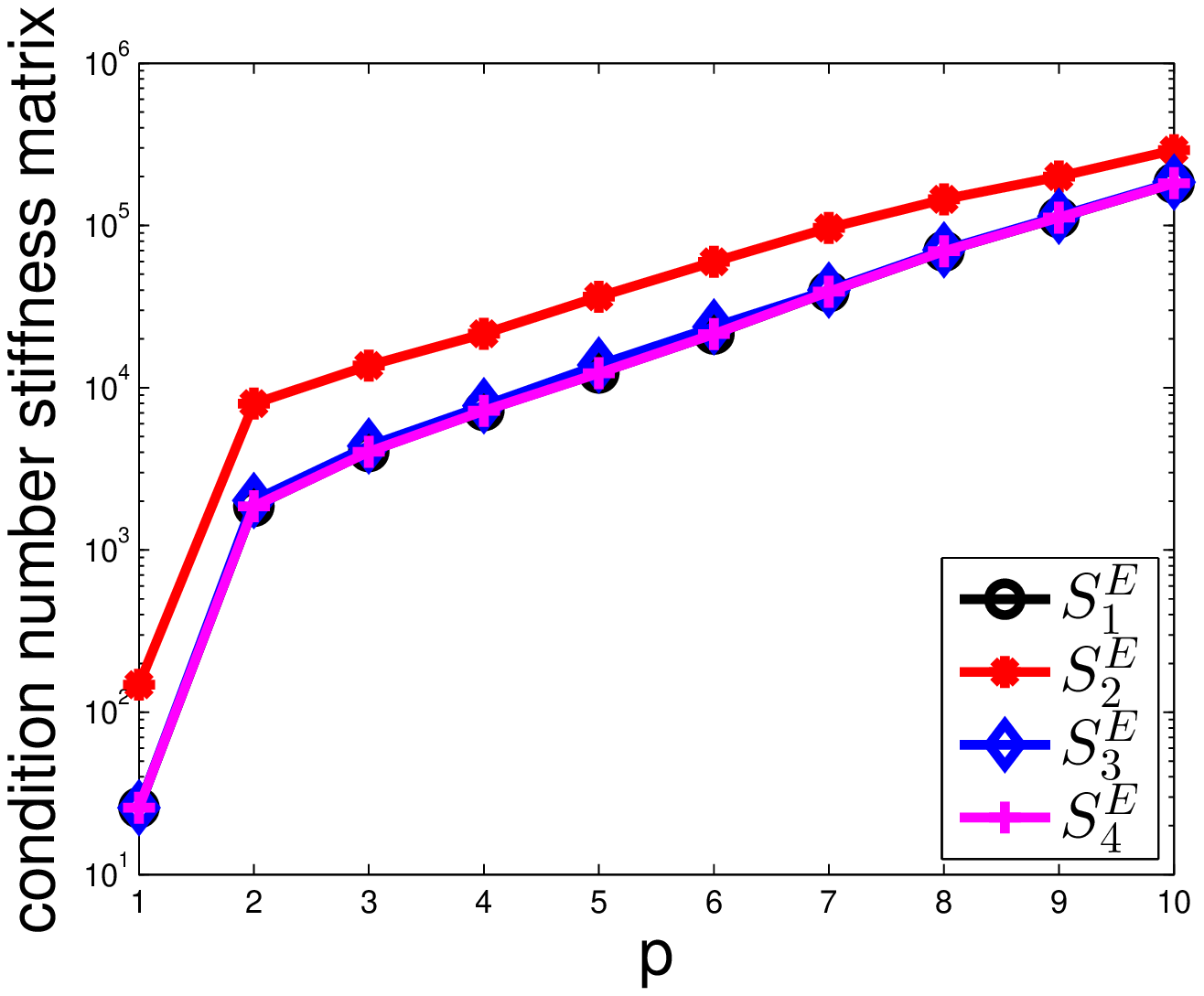}}
\subfigure {\includegraphics [angle=0, width=0.32\textwidth]{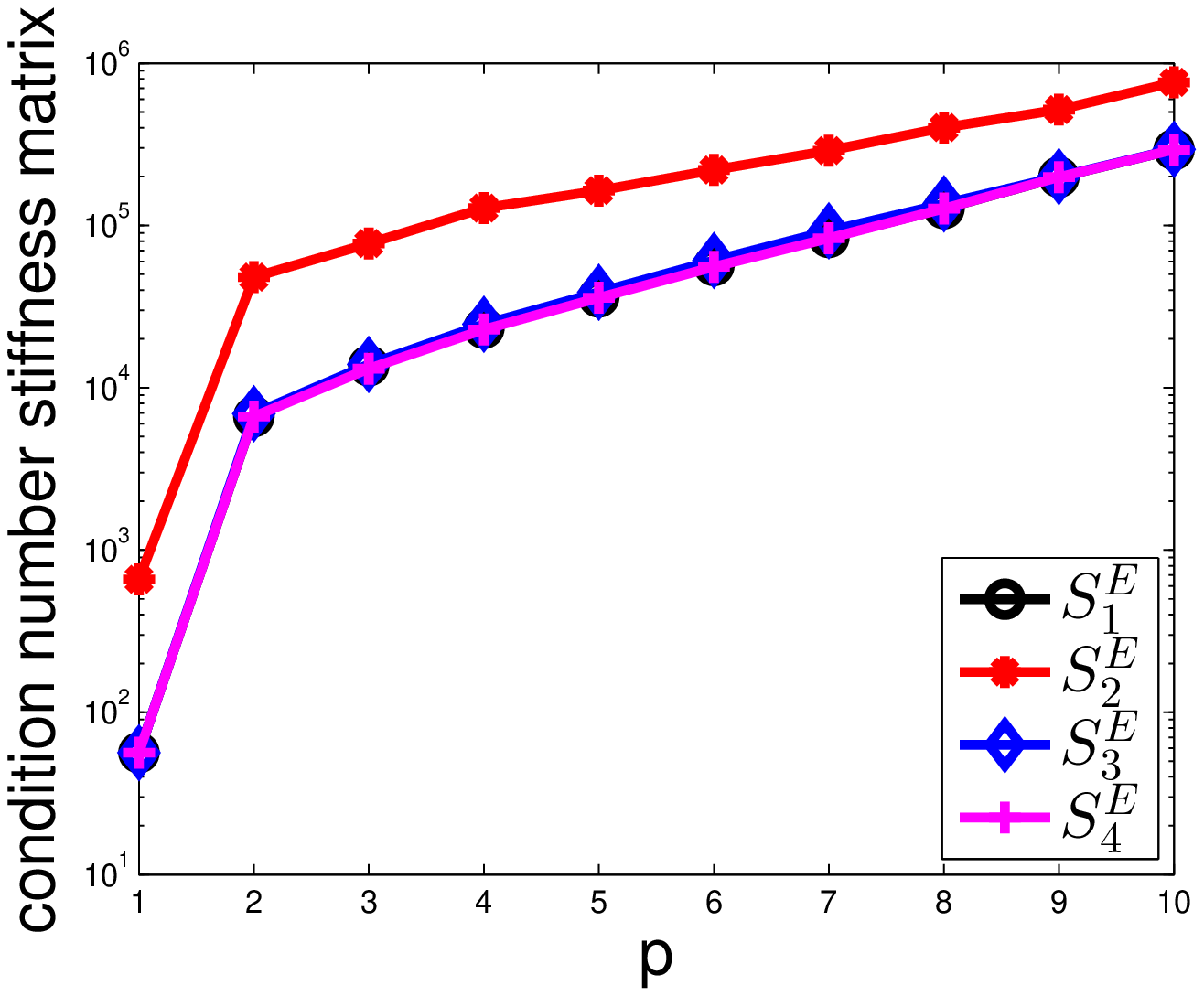}}
\caption{Condition numbers of the VEM stiffness matrix in terms of $\p$ on the three meshes depicted in Figure \ref{figure three meshes}.
The polynomial basis dual to internal moments \eqref{internal moments} is $\qalphabasisd \azpmd$. We compare the behaviour in terms of the stabilizations presented in Section \ref{subsection choices for the stabilization}.
Left: square mesh. Center: Voronoi-Lloyd mesh. Right: regular-hexagonal mesh.} \label{figure p cond various stab}
\end{figure}

In order to conclude this section, we present in Figure \ref{figure p cond various stab} a numerical test where we fix the polynomial basis dual to the internal moments \eqref{internal moments} to be $\qalphabasisd \azpmd$ and we consider the four different stabilizations
of Section \ref{subsection choices for the stabilization}.
The stabilizations, as in the case of ``bad'' geometrical deformation presented in Section \ref{subsection numerical results: collapsing polygon}, have almost the same impact on the condition number
(stabilization $\SE_2$ seems to perform slightly worse than the other stabilizations).

\subsection{Numerical results: collapsing polygons} \label{subsection numerical results: collapsing polygon}
It is also interesting to understand which is the impact of the choice of the stabilization and of the polynomial basis dual to internal moments \eqref{internal moments} in presence of
a sequence of ``bad shaped'' polygons (i.e. with collapsing bulk) on the condition number of the local stiffness matrix.
In this way, we also test the robustness of the method when assumption (\textbf{D1}) is not valid.

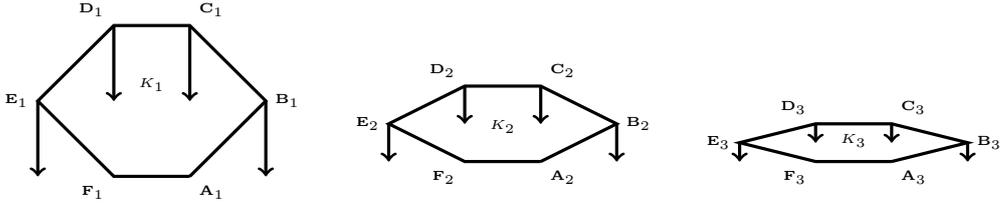
\begin{figure}  [h]
\begin{center}
\begin{minipage}{0.30\textwidth}
\begin{tikzpicture}[anchor = south]
\draw[black, very thick, -] (1,0) -- (2,1) -- (1,2) -- (0,2) -- (-1,1) -- (0,0) -- (1,0);
\draw[black, -] (1,0) node[black, below right] {\tiny{$\mathbf A_1$}}; \draw[black, -] (2,1) node[black, right] {\tiny{$\mathbf B_1$}}; \draw[black, -] (1, 2) node[black, above right] {\tiny{$\mathbf C_1$}};
\draw[black, -] (0,2) node[black, above left] {\tiny{$\mathbf D_1$}}; \draw[black, -] (-1,1) node[black, left] {\tiny{$\mathbf E_1$}}; \draw[black, -] (0,0) node[black, below left] {\tiny{$\mathbf F_1$}};
\draw[black, -] (1/2,1) node[black] {\tiny{$\E_1$}};
\draw[black, very thick, ->] (2,1) -- (2,0); \draw[black, very thick, ->] (1, 2) -- (1,1); \draw[black, very thick, ->] (0, 2) -- (0, 1); \draw[black, very thick, ->] (-1, 1) -- (-1,0);
\end{tikzpicture}
\end{minipage}
\begin{minipage}{0.30\textwidth}
\begin{tikzpicture}[anchor=south]
\draw[black, very thick, -] (1,0) -- (2, 1/2) -- (1, 1) -- (0, 1) -- (-1, 1/2) -- (0,0) -- (1,0);
\draw[white, very thick, -] (1,2) -- (0, 2);
\draw[black, -] (1, 0) node[black, below right] {\tiny{$\mathbf A_2$}}; \draw[black, -] (2,1*1/2) node[black, right] {\tiny{$\mathbf B_2$}}; \draw[black, -] (1, 2*1/2) node[black, above right] {\tiny{$\mathbf C_2$}};
\draw[black, -] (0, 2*1/2) node[black, above left] {\tiny{$\mathbf D_2$}}; \draw[black, -] (-1,1*1/2) node[black, left] {\tiny{$\mathbf E_2$}}; \draw[black, -] (0,0) node[black, below left] {\tiny{$\mathbf F_2$}};
\draw[black, -] (1/2,1*1/4) node[black] {\tiny{$\E_2$}};
\draw[black, very thick, ->] (2,1*1/2) -- (2,0); \draw[black, very thick, ->] (1, 2*1/2) -- (1,1*1/2); \draw[black, very thick, ->] (0, 2*1/2) -- (0, 1*1/2); \draw[black, very thick, ->] (-1, 1*1/2) -- (-1,0);
\end{tikzpicture}
\end{minipage}
\begin{minipage}{0.33\textwidth}
\begin{tikzpicture}[anchor=south]*
\draw[black, very thick, -] (1,0) -- (2,1*1/4) -- (1,2*1/4) -- (0,2*1/4) -- (-1,1*1/4) -- (0,0) -- (1,0);
\draw[white, very thick, -] (1,2) -- (0, 2);
\draw[black, -] (1,0) node[black, below right] {\tiny{$\mathbf A_3$}}; \draw[black, -] (2,1*1/4) node[black, right] {\tiny{$\mathbf B_3$}}; \draw[black, -] (1, 2*1/4) node[black, above right] {\tiny{$\mathbf C_3$}};
\draw[black, -] (0,2*1/4) node[black, above left] {\tiny{$\mathbf D_3$}}; \draw[black, -] (-1,1*1/4) node[black, left] {\tiny{$\mathbf E_3$}}; \draw[black, -] (0,0) node[black, below left] {\tiny{$\mathbf F_3$}};
\draw[black, -] (1/2,1*1/16) node[black] {\tiny{$\E_3$}};
\draw[black, very thick, ->] (2,1*1/4) -- (2,0); \draw[black, very thick, ->] (1, 2*1/4) -- (1,1*1/4); \draw[black, very thick, ->] (0, 2*1/4) -- (0, 1*1/4); \draw[black, very thick, ->] (-1, 1*1/4) -- (-1,0);
\end{tikzpicture}
\end{minipage}
\end{center}
\caption{First three elements of sequence $\{ \En \}_{i \in \mathbb N}$ of hexagons with collapsing bulk.}
\label{figure collapsing pentagons}
\end{figure}

For the purpose, we here present a quite limited and preliminary study. More precisely,
we consider $\{ \En \}_{i \in \mathbb N}$, sequence of ``collapsing'' hexagons, as those depicted in Figure \ref{figure collapsing pentagons}.
In particular, the coordinates of $\E_i$, the $i$-th element, are:
\begin{equation} \label{coordinates pentagon}
\mathbf A_i = (1,0),\quad \mathbf B_i = (2, 2^{-i+1}), \quad \mathbf C_i = (1, 2^{-i+2}),\quad \mathbf D _i = (0,  2^{-i+1}),\quad \mathbf E_i = (-1, 2^{-i+1}) ,\quad \mathbf F_i = (0,0).
\end{equation}
Needless to say, sequence $\En$ does not satisfy the star-shapedeness assumption (\textbf{D1}). 

In Figure \ref{figure h cond p=3,6 stab standard}, we depict the behaviour of the condition  number of the local stiffness matrix in terms of $i$, parameter used in the definition of the coordinates \eqref{coordinates pentagon} of the pentagons $\En$.
In particular, we compare such behaviour employing the three bases $\qalphabasisi \azpmd$, $i=1,2,3$, discussed in Section \ref{subsection choices for the basis} and choosing $\p=3$ and $\p=6$, respectively.
The stabilization is fixed to be $\SE_1$ defined in \eqref{classical stabilization}.
\begin{figure}  [h]
\centering
\subfigure {\includegraphics [angle=0, width=0.45\textwidth]{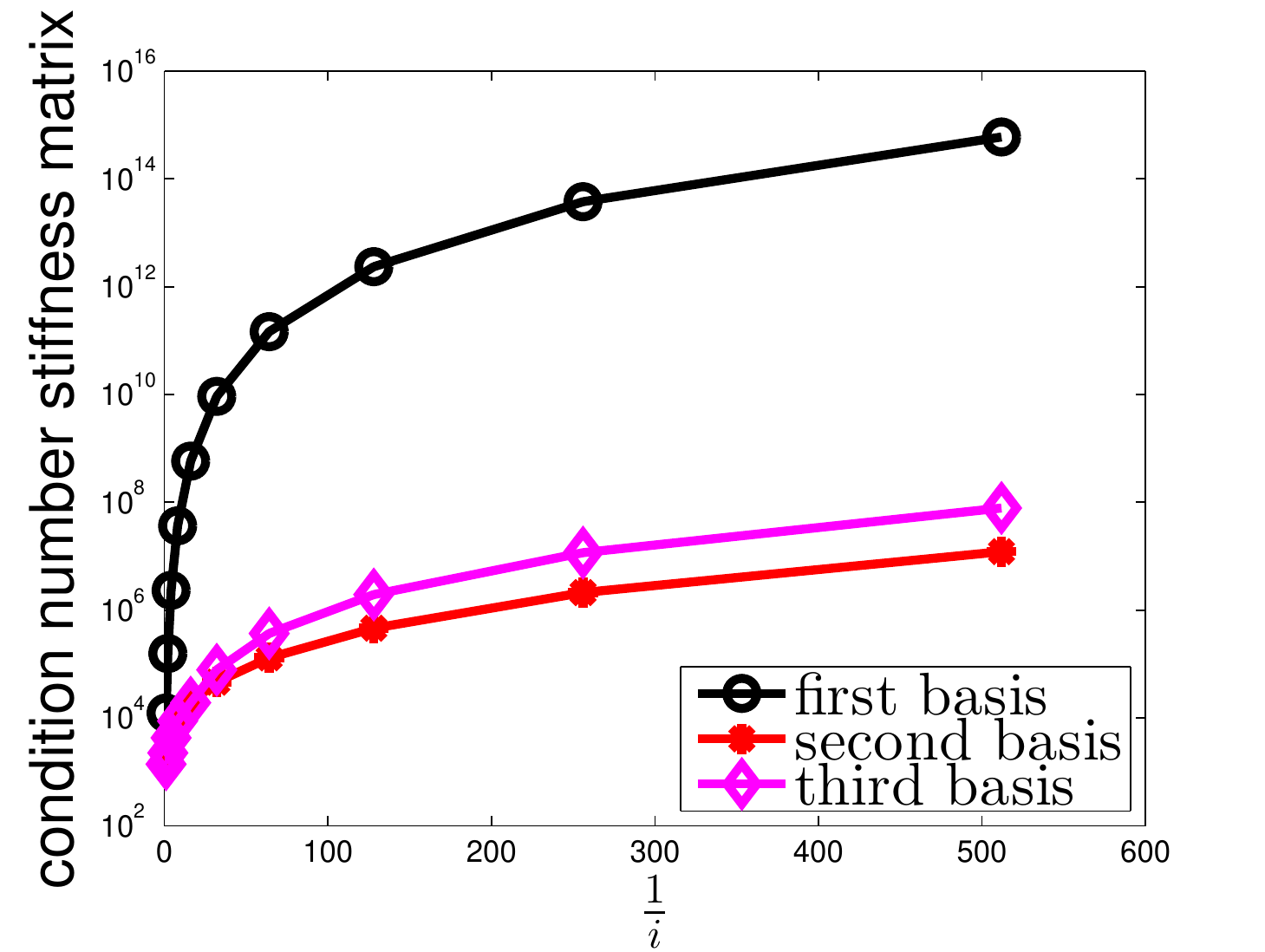}}
\subfigure {\includegraphics [angle=0, width=0.45\textwidth]{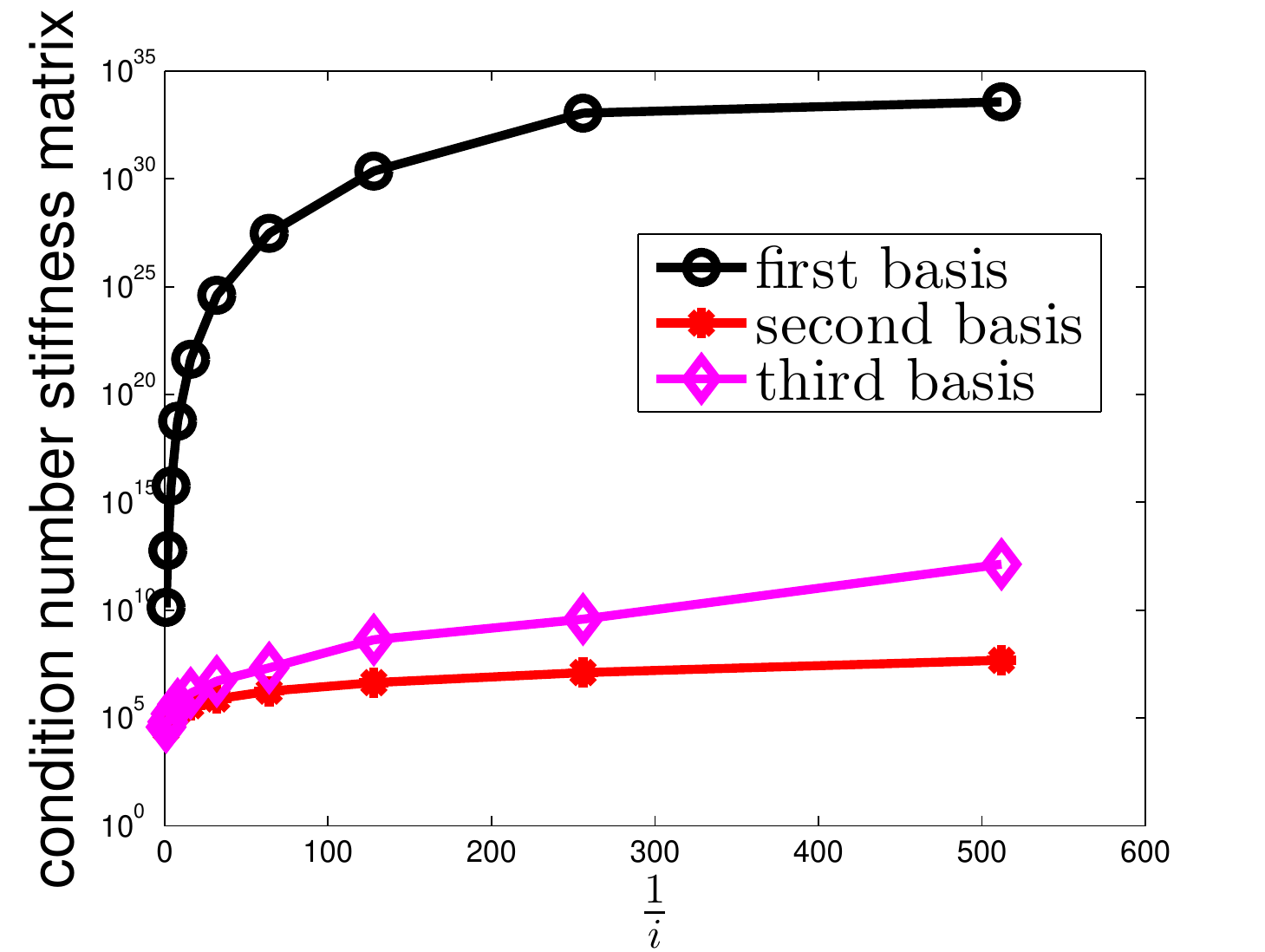}}
\caption{Condition numbers of the local VEM stiffness matrix on the sequence of hexagons depicted in Figure \ref{figure collapsing pentagons} in dependence of $i$, parameter used in the definition of the coordinates \eqref{coordinates pentagon} of the pentagons $\En$.
The stabilization is fixed and equal to $\SE_1$ \eqref{classical stabilization}. We compare the behaviour in terms of the three polynomial bases presented in Section \ref{subsection choices for the basis}.
Left $\p$=3. Right: $\p$=6.} \label{figure h cond p=3,6 stab standard}
\end{figure}
From Figure \ref{figure h cond p=3,6 stab standard}, we deduce that the standard choice for the polynomial basis \eqref{first choice basis} leads to a dramatic growth of the condition number.
It turns out that the safest choice, in terms of ill-conditioning, is the one associated with basis $\qalphabasisd \azpmd$, which we recall is obtained by an orthonormalization of the standard monomial basis $\qalphabasisu \azpmd$ via a stable Gram-Schmidt process.
Basis $\qalphabasist \azpmd$ , although behaves much better than the monomial basis, is not as good as $\qalphabasisd \azpmd$.

\begin{figure}  [h]
\centering
\subfigure {\includegraphics [angle=0, width=0.45\textwidth]{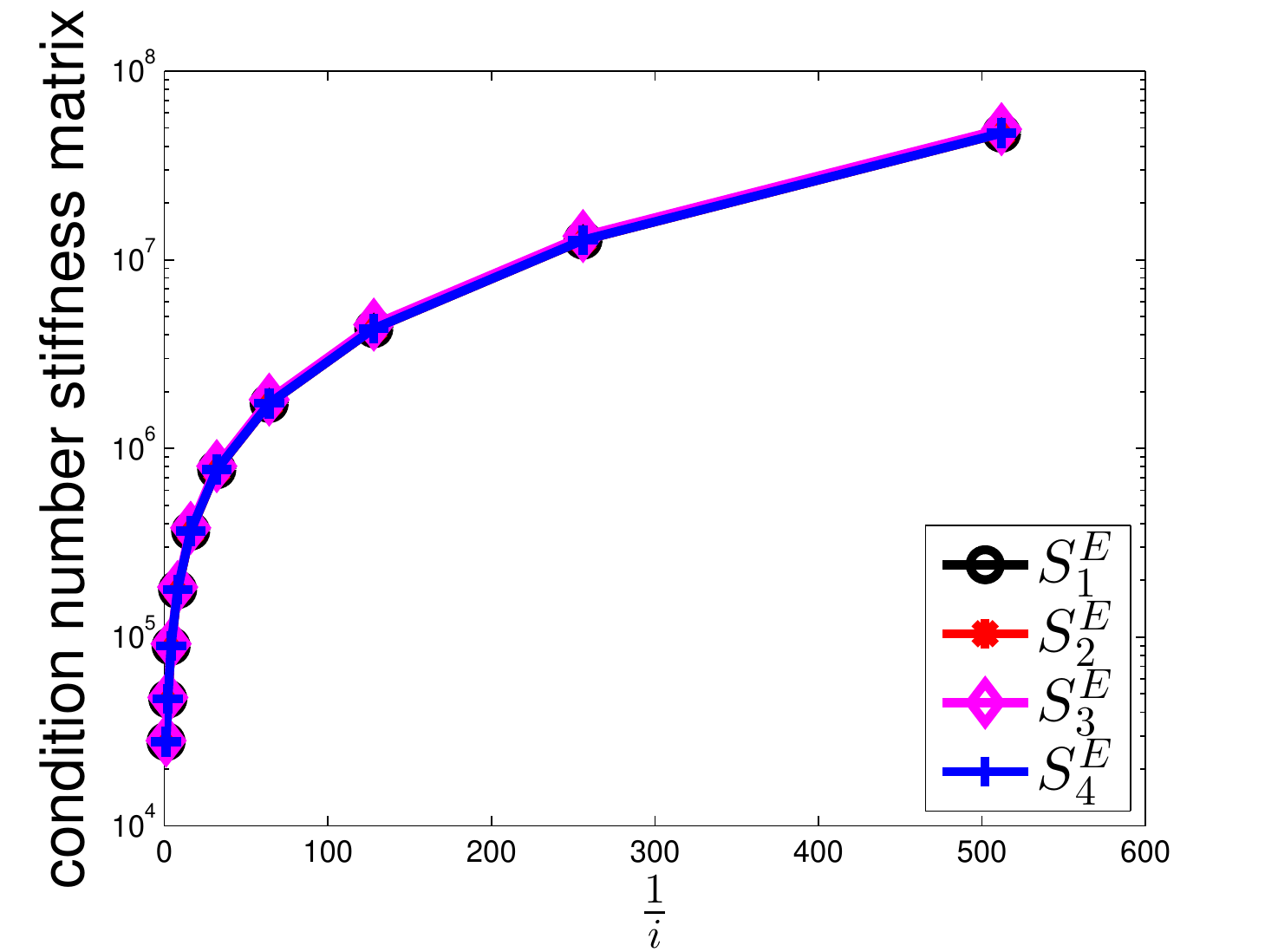}}
\caption{Condition numbers of the local VEM stiffness matrix on the sequence of hexagons depicted in Figure \ref{figure collapsing pentagons} in dependence of $i$,
parameter used in the definition of the coordinates \eqref{coordinates pentagon} of sequence $\{\En\}_{i\in \mathbb N}$. The polynomial basis, dual to the internal moments \eqref{internal moments} is fixed to be $\qalphabasisd \azpmd$. 
We compare the behaviour in terms of the four stabilizations  presented in Section \ref{subsection choices for the stabilization}. The degree of accuracy is $\p$=6.} \label{figure h cond p=6 various stab bases AMV}
\end{figure}

Next, in Figure \ref{figure h cond p=6 various stab bases AMV}, we compare the condition number of the stiffness matrix by fixing $\p=6$ and the polynomial basis $\qalphabasisd \azpmd$, which, from the previous tests, seems to be the best for the conditioning of VEM,
and by modifying the choice of the stabilizations; more precisely, we will consider the four stabilization discussed in Section \ref{subsection choices for the stabilization}.
We deduce from Figure \ref{figure h cond p=6 various stab bases AMV} that the choice of the stabilization does not have evident effects on the condition number, at least for the presented tests.

\begin{figure}  [h]
\centering
\subfigure {\includegraphics [angle=0, width=0.45\textwidth]{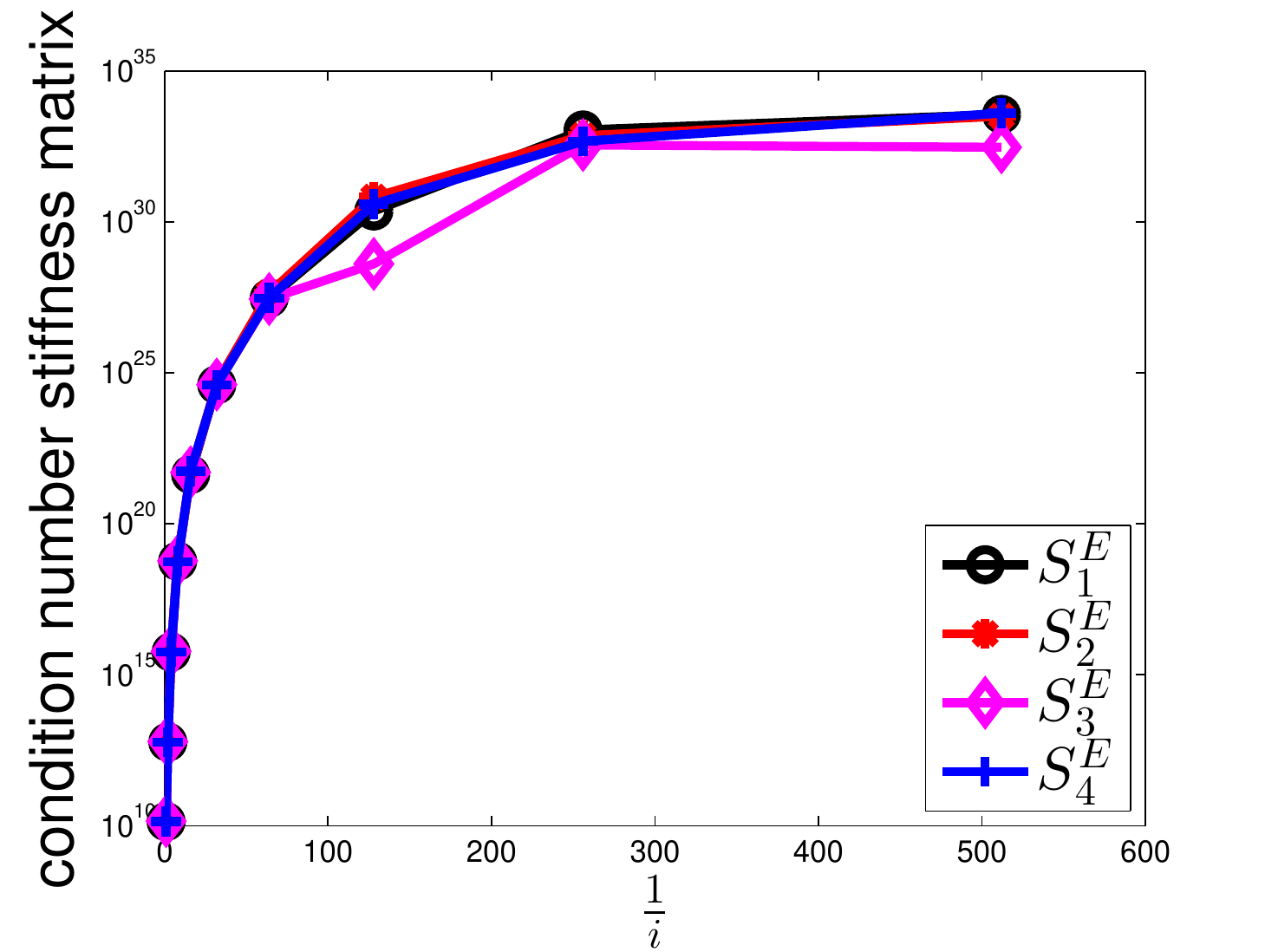}}
\caption{Condition numbers of the local VEM stiffness matrix on the sequence of hexagons depicted in Figure \ref{figure collapsing pentagons} in dependence of $i$,
parameter used in the definition of the coordinates \eqref{coordinates pentagon} of sequence $\{\En\}_{i\in \mathbb N}$. The polynomial basis, dual to the internal moments \eqref{internal moments} is fixed to be $\qalphabasisu \azpmd$. 
We compare the behaviour in terms of the four stabilizations  presented in Section \ref{subsection choices for the stabilization}. The degree of accuracy is $\p$=6.} \label{figure h cond p=6 various stab standard basis}
\end{figure}

As a byproduct, in Figure \ref{figure h cond p=6 various stab standard basis} we consider a comparison between the four stabilizations by employing the standard monomial basis $\qalphabasisu \azpmd$ as dual basis for internal moments \eqref{internal moments}.
Again, we assume $\p=6$.

\subsection{Numerical results: hanging nodes} \label{subsection numerical results: hanging nodes}
As a final set of numerical results, we study the behaviour of the condition number of the stiffness matrix employing various bases and stabilizations in presence of hanging nodes collapsing on a vertex,
checking thus the robustness of the method when assumption (\textbf{D2}) is not fulfilled.

\begin{figure}  [h]
\begin{center}
\begin{minipage}{0.30\textwidth}
\begin{tikzpicture} [scale=3]
\draw[black, very thick, -] (1,0) -- (1,1) -- (0.5,1) -- (0,1) -- (0,0) -- (1,0);
\draw[black, -] (1,0) node[black, below right] {\tiny{$\mathbf A_1$}}; \draw[black, -] (1,1) node[black, above right] {\tiny{$\mathbf B_1$}}; \draw[black, -] (0.5, 1) node[black, above] {\tiny{$\mathbf C_1$}};
\draw[black, -] (0,1) node[black, above left] {\tiny{$\mathbf D_1$}}; \draw[black, -] (0,0) node[black, below left] {\tiny{$\mathbf E_1$}};
\draw[black, very thick, ->] (1/2,1.05) -- (1/4,1.05);
\draw[black, fill=black] (0.5,1) circle (0.3mm);
\end{tikzpicture}
\end{minipage}
\begin{minipage}{0.30\textwidth}
\begin{tikzpicture} [scale=3]
\draw[black, very thick, -] (1,0) -- (1,1) -- (1/4,1) -- (0,1) -- (0,0) -- (1,0);
\draw[black, -] (1,0) node[black, below right] {\tiny{$\mathbf A_2$}}; \draw[black, -] (1,1) node[black, above right] {\tiny{$\mathbf B_2$}}; \draw[black, -] (1/4, 1) node[black, above] {\tiny{$\mathbf C_2$}};
\draw[black, -] (0,1) node[black, above left] {\tiny{$\mathbf D_2$}}; \draw[black, -] (0,0) node[black, below left] {\tiny{$\mathbf E_2$}};
\draw[black, very thick, ->] (1/4,1.05) -- (1/8,1.05);
\draw[black, fill=black] (1/4,1) circle (0.3mm);
\end{tikzpicture}
\end{minipage}
\begin{minipage}{0.3\textwidth}
\begin{tikzpicture} [scale=3]
\draw[black, very thick, -] (1,0) -- (1,1) -- (1/8,1) -- (0,1) -- (0,0) -- (1,0);
\draw[black, -] (1,0) node[black, below right] {\tiny{$\mathbf A_3$}}; \draw[black, -] (1,1) node[black, above right] {\tiny{$\mathbf B_3$}}; \draw[black, -] (1/8, 1) node[black, above right] {\tiny{$\mathbf C_3$}};
\draw[black, -] (0,1) node[black, above left] {\tiny{$\mathbf D_3$}}; \draw[black, -] (0,0) node[black, below left] {\tiny{$\mathbf E_3$}};
\draw[black, very thick, ->] (1/8,1.05) -- (1/16,1.05);
\draw[black, fill=black] (1/8,1) circle (0.3mm);
\end{tikzpicture}
\end{minipage}
\end{center}
\caption{First three elements of sequence $\{ \En \}_{i \in \mathbb N}$ with hanging node collapsing on a vertex.}
\label{figure square pentagons}
\end{figure}
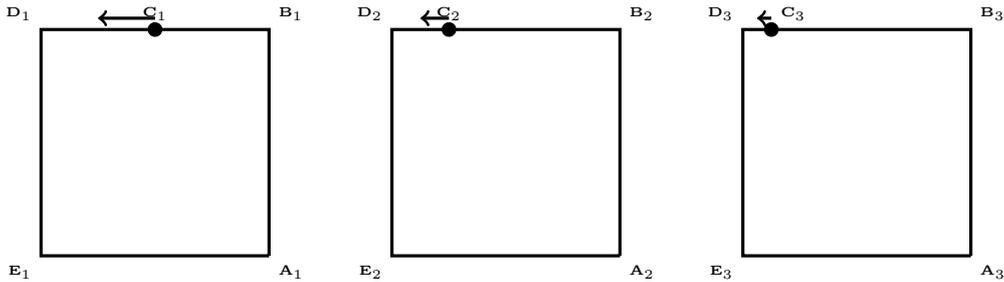

Again, we present here only a quite limited and preliminary study.
In particular, we present a sequence of ``squared pentagons'', that is a sequence of squares with a hanging node on a prescribed edge.
More precisely, see Figure \ref{figure square pentagons}, we consider a sequence $\{\E_i\}_{i\in \mathbb N}$ such that each $\E_i$, $i\in \mathbb N$, has the following set of coordinates:

\begin{equation} \label{coordinates squared pentagon}
\mathbf A_i = (1,0),\quad \mathbf B_i = (1,1), \quad \mathbf C_i = (2^{-i}, 1),\quad \mathbf D _i = (0, 1),\quad \mathbf E_i = (0,0) ,\, i\in \mathbb N.
\end{equation}

In Figure \ref{figure h cond p=3,6 stab standard hanging nodes}, we depict the behaviour of the condition  number of the local stiffness matrix in terms of $i$,
parameter used in the definition of the coordinates \eqref{coordinates squared pentagon} of the pentagons $\En$.
In particular, we compare such behaviour employing the three bases $\qalphabasisi \azpmd$, $i=1,2,3$, discussed in Section \ref{subsection choices for the basis} and choosing $\p=3$ and $\p=6$, respectively.
The stabilization is fixed to be $\SE_1$ defined in \eqref{classical stabilization}.
\begin{figure}  [h]
\centering
\subfigure {\includegraphics [angle=0, width=0.45\textwidth]{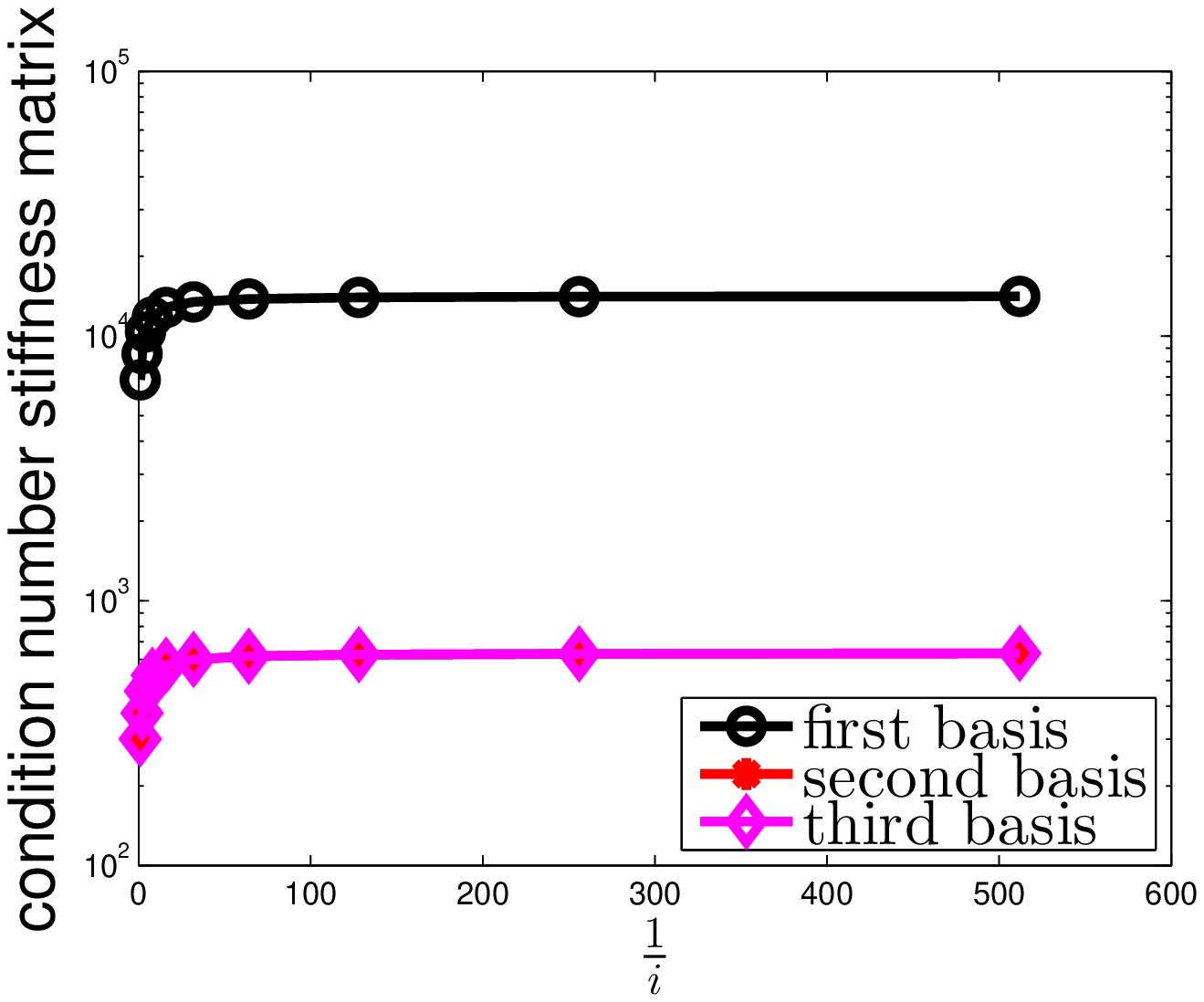}}
\subfigure {\includegraphics [angle=0, width=0.45\textwidth]{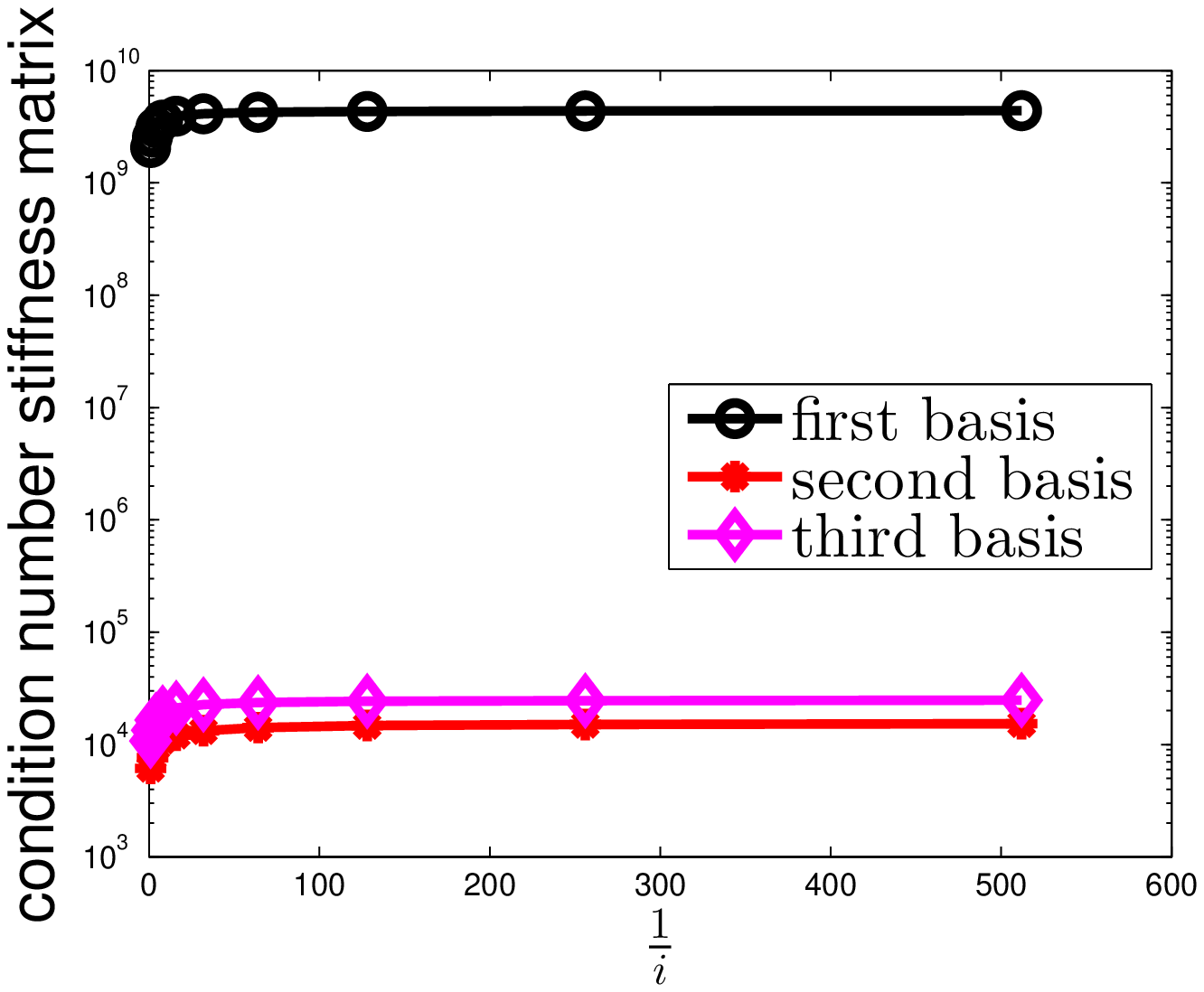}}
\caption{Condition numbers of the local VEM stiffness matrix on the sequence of pentagons (squares with a hanging node)
depicted in Figure \ref{figure square pentagons} in dependence of $i$, parameter used in the definition of the coordinates \eqref{coordinates pentagon} of the pentagons $\En$.
The stabilization is fixed and equal to $\SE_1$ \eqref{classical stabilization}. We compare the behaviour in terms of the three polynomial bases presented in Section \ref{subsection choices for the basis}.
Left $\p$=3. Right: $\p$=6. \label{figure h cond p=3,6 stab standard hanging nodes}}
\end{figure}
The condition number is almost independent of parameter $i$ for all choices of the canonical basis. This is not surprising
since the bulk of the elements in the sequence remains the same for all $i$.
However, when employing basis $\qalphabasisu \azpmd$, the condition number is higher.

\begin{figure}  [h]
\centering
\subfigure {\includegraphics [angle=0, width=0.45\textwidth]{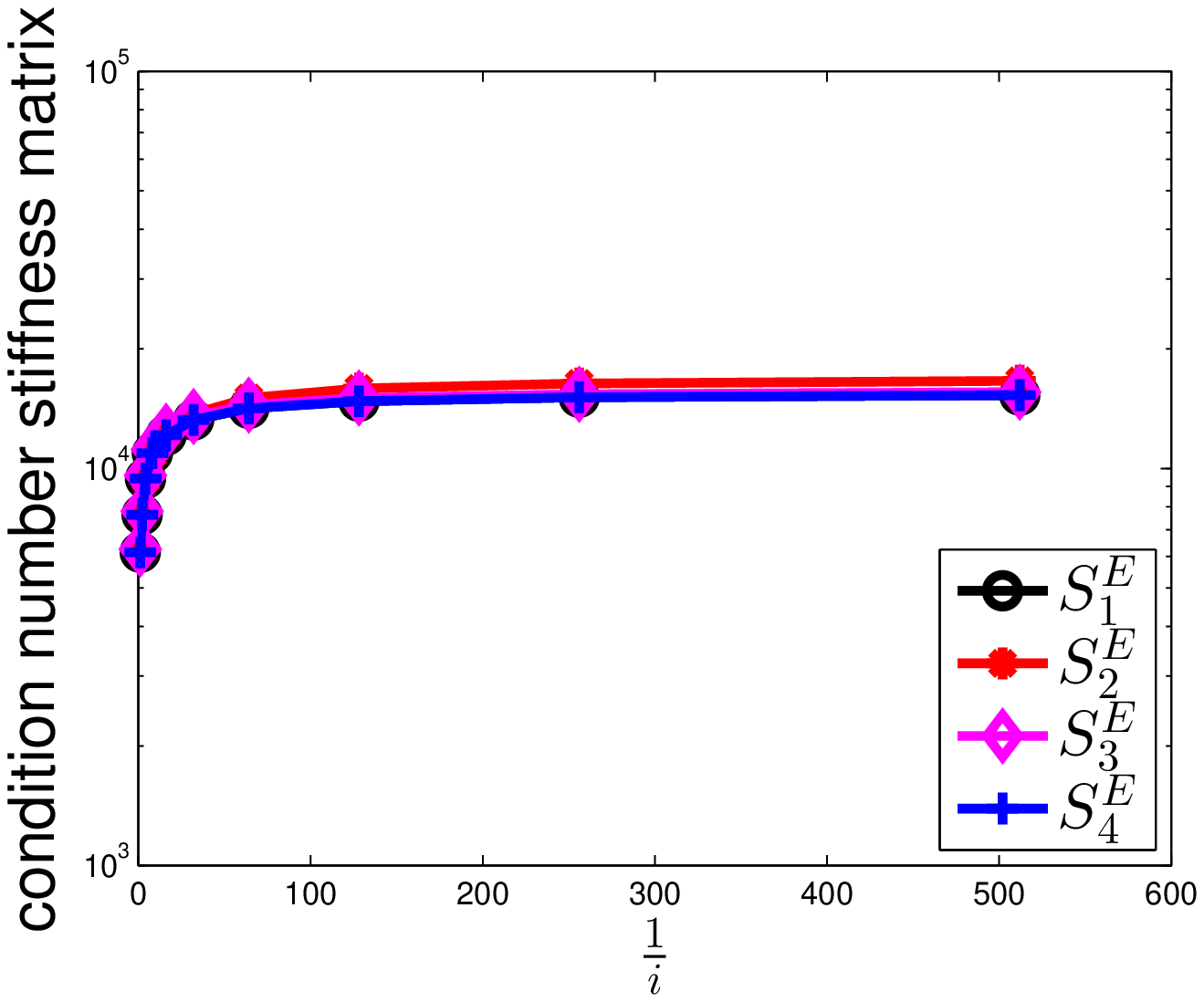}}
\caption{Condition numbers of the local VEM stiffness matrix on the sequence of pentagons (squares with a hanging node)
depicted in Figure \ref{figure square pentagons} in dependence of $i$,
parameter used in the definition of the coordinates \eqref{coordinates pentagon} of sequence $\{\En\}_{i\in \mathbb N}$. The polynomial basis, dual to the internal moments \eqref{internal moments} is fixed to be $\qalphabasisd \azpmd$. 
We compare the behaviour in terms of the four stabilizations  presented in Section \ref{subsection choices for the stabilization}. The degree of accuracy is $\p$=6. \label{figure h cond p=6 various stab bases AMV hanging nodes}}
\end{figure}

\begin{figure}  [h]
\centering
\subfigure {\includegraphics [angle=0, width=0.45\textwidth]{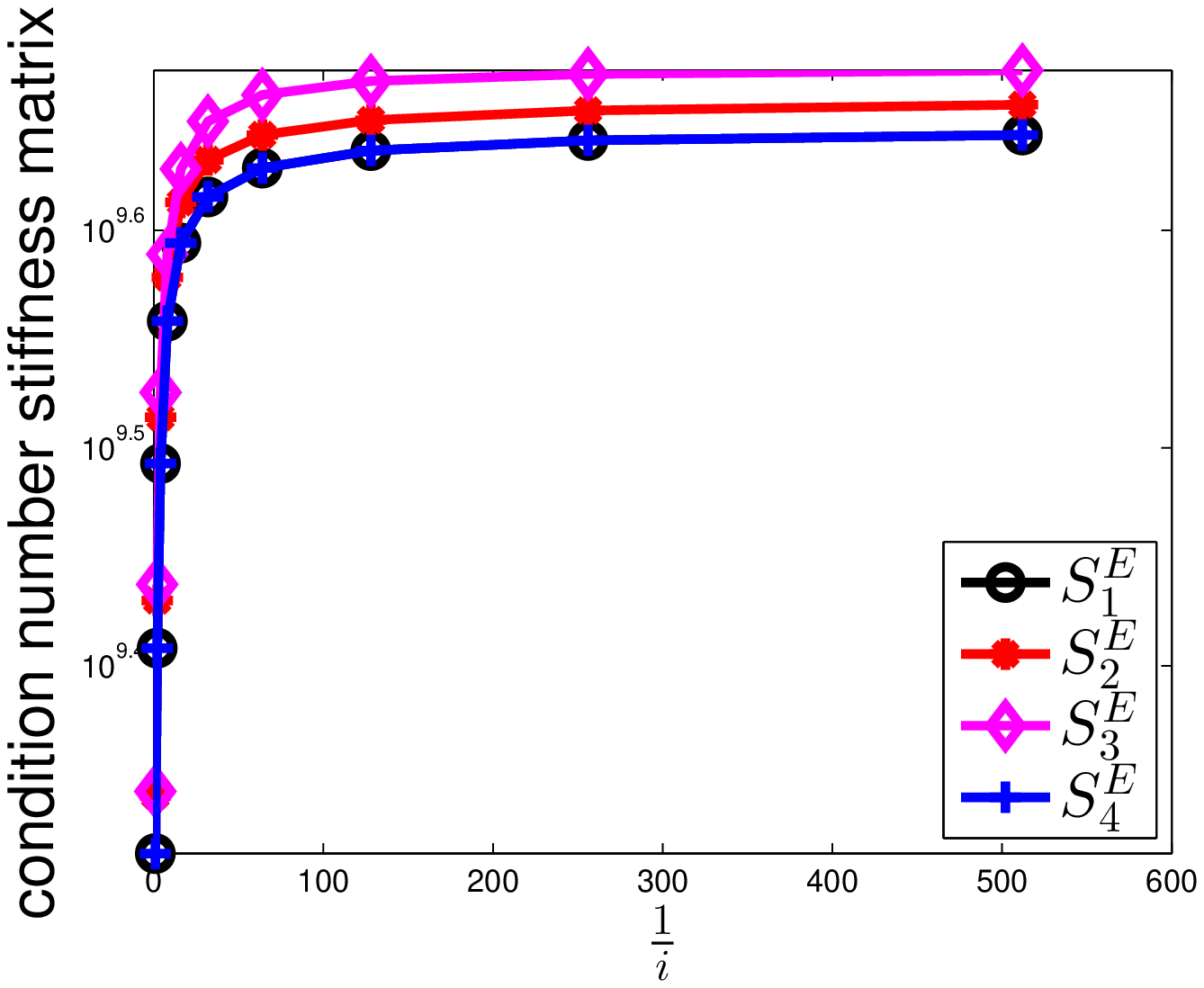}}
\caption{Condition numbers of the local VEM stiffness matrix on the sequence of pentagons (squares with a hanging node)
depicted in Figure \ref{figure square pentagons} in dependence of $i$,
parameter used in the definition of the coordinates \eqref{coordinates pentagon} of sequence $\{\En\}_{i\in \mathbb N}$. The polynomial basis, dual to the internal moments \eqref{internal moments} is fixed to be $\qalphabasisu \azpmd$. 
We compare the behaviour in terms of the four stabilizations  presented in Section \ref{subsection choices for the stabilization}. The degree of accuracy is $\p$=6. \label{figure h cond p=6 various stab standard basis hanging nodes}}
\end{figure}

Mimicking what done in Section \ref{subsection numerical results: collapsing polygon}, we compare in Figures \ref{figure h cond p=6 various stab bases AMV hanging nodes} and \ref{figure h cond p=6 various stab standard basis hanging nodes}
the condition number of the stiffness matrix by fixing $\p=6$ and the polynomial bases $\qalphabasisd \azpmd$ and $\qalphabasisu \azpmd$, respectively, 
and by considering the four stabilization discussed in Section \ref{subsection choices for the stabilization}.

We deduce from Figures \ref{figure h cond p=6 various stab bases AMV hanging nodes} and \ref{figure h cond p=6 various stab standard basis hanging nodes} again that the behaviour of the method employing different stabilizations
is basically the same.

\section {Conclusions} \label{section conclusions}
In the present work, we addressed and suggested possible cures to the problem of the ill-conditioning of the Virtual Element Method,
arising from high values of the polynomial degree $\p$ and in presence of highly anisotropic elements.
In particular, we focused our attention on the effects of the stabilization of the method and the choice of internal degrees of freedom.
It turned out that, whereas various stabilizations presented in literature have almost the same effect on the condition number of the stiffness matrix, 
the choice of the internal degrees of freedom has a deep impact.

We suggested two practical modifications of such internal degrees of freedom which greatly improve the behaviour of high-order VEM and VEM in presence of bad-shaped polygons.

It is worth to mention that we focused our attention to a simple 2D Poisson problem only.
If one turns for instance to 3D problems, then the choice of stabilization \eqref{stability S} plays a major role, as shown in \cite{preprint_VEM3Dbasic, fetishVEM3D}.

\appendix
\section{Appendix: a hitchhikers guide for the implementation of VEM using basis $\qalphabasisd \azpmd$ } \label{section hitchhikers}
We discuss here the implementation details for the construction of the method employing basis $\qalphabasisd \azpmd$ as a dual basis for the internal moments \eqref{internal moments}.
Henceforth, we adopt the notation of Section \ref{subsection choices for the basis}.
Moreover, given a matrix $\A \in \mathbb R^{n\times m}$, we will denote by:
\[
\A(i:j, \ell: k ),
\]
the submatrix of $\A$ from row $i$ to $j$ and from column $\ell$ to $k$.

Let us define by $\qalphabasisd \azp$ the basis of $\mathbb P_\p(\E)$ obtained by an $L^2$ ortonormalization of basis $\qalphabasisu \azpmd$ using the stable Gram-Schimdt process presented in \cite{BassiBottiColomboDipietroTesini}.
In particular, we can write (always using with a little abuse of notation the bijection \eqref{natural bijection}):
\[
\qalphad(\xbold) = \sum_{\beta = 1}^{\alpha} \GS_{\alpha,\beta}\, \qalphau (\xbold) \quad \quad \forall \boldalpha \in \mathbb N^2,\quad \vert \boldalpha \vert = 0, \dots, \p,
\]
where $\GS$ is the lower triangular matrix containing the orthonormalizing coefficients.

In \cite{hitchhikersguideVEM}, the implementation details employing basis $\qalphabasisu \azpmd$ defined in \eqref{first choice basis} were discussed. In particular, it was proven that the local stiffness matrix was defined through some auxiliary matrices which we recall here:
\begin{equation} \label{VEM matrices part1}
\begin{split}
&\G_{\alpha,\beta} = 	\begin{cases}
					P_0(\qbetau) & \text{if } \alpha =1\\
					(\nabla \qalphau, \nabla \qbetau)_{0,\E} & \text{if } \alpha \ge 2\\
				\end{cases}, \quad \quad
\quad \Gtilde _{\alpha, \beta} = (\nabla \qalphau, \nabla \qbetau)_{0,\E},\\
&\D_{j,\alpha} = \dof_j^1 (\qalphau), \quad \quad \quad 	\B_{\alpha,j} = \begin{cases}
															P_0 \varphi^1_j & \text{if } \alpha=1\\
															(\nabla \qalphau, \nabla \varphi_j^1) _{0,\E} & \text{if } \alpha \ge 2\\
														\end{cases},\\
&\forall \boldalpha,\, \boldbeta \in \mathbb N^2,\; \vert \boldalpha \vert,\, \vert \boldbeta \vert = 0,\dots, \p,\; i=0,\dots,\NdofE,
\end{split}
\end{equation}
where $P_0$ is defined in \eqref{fixing constants}.

The local stiffness matrix reads:
\begin{equation} \label{local stiffness matrix}
\K_\p^\E = (\Pinablastar)^t \cdot \Gtilde \cdot (\Pinablastar) + (\Id - \Pinablamat)^t \cdot \SEmat \cdot (\Id - \Pinablamat),
\end{equation}
where $\SEmat$ denotes the matrix associated with any of the bilinear forms $\SE$ introduced in \eqref{stability S}, where $\Pinablastar$ denotes the matrix associated with operator $\Pinablap$ introduced in \eqref{H1 orthogonality}
acting from the local VE space $\VnE$ to $\mathbb P_\p(\E)$ with respect to basis $\qalphabasisu \azpmd$ and where $\Pinablamat$ denotes the matrix associated with the operator $\Pinablap$ introduced in \eqref{H1 orthogonality}
acting from the local VE space $\VnE$ to $\mathbb P_\p(\E)$ with respect to basis $\{\varphi_j^1\}_{j=1}^{\NdofE}$.

It was shown in \cite{hitchhikersguideVEM} that:
\[
\Pinablastar = \G^{-1} \cdot \B,\quad \quad \quad \Pinablamat = \D \cdot \G^{-1} \cdot \B.
\]
The aim of the present section is to write the local stiffness matrix employing the new canonical basis associated with polynomial basis $\qalphabasisd \azpmd$
and expanding all the projectors with respect to the polynomial spaces on polynomial basis $\qalphabasisd \azp$.

For the sake of simplicity, we denote the counterpart of the VEM matrices in \eqref{VEM matrices part1} and \eqref{local stiffness matrix} associated with bases $\qalphabasisd \azp$ and $\{\varphi_j^2\}_{j=1}^{\NdofE}$, with a bar at the top of each one of them.

We explain how to compute in the new setting such new matrices.
We start with matrix $\Gtildebar$, which is defined as:
\[
\Gtildebar_{\alpha, \beta} = (\nabla \qalphad, \nabla \qbetad)_{0,\E}\quad \quad \forall \boldalpha,\, \boldbeta \in \mathbb N^2,\quad \vert \boldalpha \vert, \, \vert \boldbeta \vert  = 0,\dots,\p.
\]
One simply has to compute:
\[
\Gtildebar = \GS \cdot \Gtilde \cdot \GS^T.
\]
Now, we consider matrix $\Gbar$ defined as:
\[
\Gbar_{\alpha, \beta} = \begin{cases}
P_0(\qbetad) & \text{if } \alpha = 1\\
(\nabla \qalphad, \nabla \qbetad) & \text{if } \alpha \ge 2\\
\end{cases}\quad \quad \forall \boldalpha,\, \boldbeta \in \mathbb N^2,\quad \vert \boldalpha \vert, \, \vert \boldbeta \vert  = 0,\dots,\p,
\]
where $P_0(\cdot)$ is defined in \eqref{fixing constants}. We obviously have to take care only of the first line of the matrix since the remainder is inherited from $\Gtildebar$.
We distinguish two cases.
\begin{itemize}
\item [($\p = 1$)] $P_0 (\qbetad) = \frac{1}{\NE} \sum_{\ell=1}^{\NE} \qbetad (\nu_\ell)$, where we recall that $\nu_\ell$ is the $\ell$-th vertex of $\E$. We shall then write:
\[
\begin{split}
P_0(\qbetad) 	& = \frac{1}{\NE} \sum _{\ell=1}^{\NE} \left(\sum_{\gamma \le \beta} \GS_{\beta, \gamma} \qgamma^1 (\nu_\ell) \right) = \sum_{\gamma =1}^{\beta} \GS_{\beta,\gamma} \left( \frac{1}{\NE} \sum_{\ell=1}^{\NE} \qgamma^1 (\nu_\ell)  \right)\\
			& = \sum_{\gamma =1}^{\beta} \GS_{\beta, \gamma} \G_{1,\gamma} = \GS(\beta, 1:\beta) \cdot \G(1,1:\beta)^T.
\end{split}
\]
\item[($\p \ge 2$)] In this case, we have:
\[
P_0(\qbetad) = \frac{1}{\vert \E \vert} \int_\E \qbetad = \GS_{1,1}^{-1} \frac{1}{\vert \E \vert} \int _\E \qbetad \GS_{1,1} = \GS_{1,1}^{-1} \frac{1}{\vert \E \vert} \int _\E \q_1^2 \qbetad  = \begin{cases}
																							\frac{1}{\GS_{1,1} \vert \E \vert} & \text{if } \beta=1\\
																							0                         & \text{else}\\
																						\end{cases},
\]
since basis $\{\qalphad\}_{\alpha=1}^{\np}$ is $L^2(\E)$ orthonormal by construction.
\end{itemize}
\medskip\medskip\medskip

Next, we turn our attention to matrix $\Dbar$ which is defined as:
\[
\Dbar_{j,\alpha} = \dof_j^2(\qalphad)\quad \quad \forall \boldalpha \in \mathbb N^2,\quad \vert \boldalpha \vert = 0,\dots, \p,\quad j=1,\dots, \NdofE
\]
and we distinguish two situations.

\begin{itemize}
\item Let us consider firstly the boundary dofs:
\[
\Dbar_{j,\alpha} = \dof^2_j(\qalphad) = \qalphad(\xi_j) = \sum_{\beta =1}^{\alpha} \GS_{\alpha, \beta} \qbetau(\xi_j) = \sum_{\beta =1}^{\alpha} \GS_{\alpha, \beta} \D_{j,\beta} = \GS (\alpha, 1: \alpha) \cdot \D(j, 1:\alpha)^T.
\]
where $\xi_j$ is a proper node on the boundary.
\item
Next, we deal with the internal dofs. One simply has, if $\qgammad$ is the polynomial associated with $\dof_j$:
\[
\dof^2_j(\qalphad) = \frac{1}{\vert \E \vert} \int_\E \qalphad \qgamma^2 = \frac{1}{\vert \E \vert} \delta_{\alpha,j}\quad \quad 	\forall j =1,\dots, \npmd,
\]
owing again to the $L^2$ orthonormality of basis $\qalphabasisd \azp$.
\end{itemize}
\medskip\medskip\medskip

Finally, we discuss the construction of matrix $\Bbar$ which is defined as:
\[
\Bbar_{\alpha,j} = 	\begin{cases}
				P_0(\varphi_j^2) & \text{if } \alpha=1\\
				(\nabla \qalphad, \nabla \varphi_j^2) & \text{if } \alpha\ge 2\\
			\end{cases}\quad \quad \forall \boldalpha \in \mathbb N^2,\quad \vert \boldalpha \vert = 0,\dots,\p,\quad \forall j = 1,\dots, \NdofE.
\]
We firstly deal with the first line and we consider the two cases $\p=1$ and $\p \ge 2$.
\begin{itemize}
\item[($\p=1$)]  $P_0(\varphi^2_j) = \frac{1}{\NE} \sum_{\ell=1}^{\NE} \varphi_j^2(\nu_\ell) = \frac{1}{\NE} \sum_{\ell=1}^{\NE} \varphi_j^1(\nu_\ell)$,
where we recall that $\{\nu_\ell\}_{\ell =1}^{\NE}$ is the set of vertices of polygon $\E$.
Thus $\Bbar_{1,j} = \B_{1,j}$ $\forall j = 1,\dots, \NdofE$.
\item[($\p\ge 2$)] In this case, we can write:
\[
\begin{split}
P_0(\varphi _j^2) 	&= \frac{1}{\vert \E \vert } \int_\E \varphi _j^2 = \GS_{1,1}^{-1} \frac{1}{\vert \E \vert} \int_\E \varphi _j^2\, \GS_{1,1} =  \GS_{1,1}^{-1} \frac{1}{\vert \E \vert} \int_\E \varphi_j^2\, \q_1^2\\
			& = 	\begin{cases}
					\GS_{1,1}^{-1} & \text{if } \varphi_j \text{ is the first internal element},\\
					0 & \text{else},\\
				\end{cases}
\end{split}
\] 
since $\q_1^2 = \GS_{1,1} \q_1^1 = \GS_{1,1}$. Thus $\Bbar_{1,j} = \GS_{1,1}^{-1} \B_{1,j}$ $\forall j = 1,\dots, \NdofE$.
\end{itemize}
Next, we treat all the other lines. We must compute $(\nabla \qalphad, \nabla \varphi_j^2)_{0,\E}$. Again, we consider two different situations.
\begin{itemize}
\item If $\varphi_j^2$ is a boundary basis function, i.e. $j=1,\dots, \p \NE$, where we recall that $\NE$ is the number of edges (and vertices) of $\E$, then:
\[
\begin{split}
(\nabla \qalphad, \nabla \varphi_j^2)_{0,\E} 		& = \int_{\partial \E} (\partial_\n \qalphad) \,\varphi_j^2 = \sum_{\beta = 1}^{\alpha} \GS_{\alpha, \beta} \int _{\partial \E} (\partial _\n \qbeta^1) \, \varphi_j^2\\
								& = \sum_{\beta =1}^{\alpha} \GS_{\alpha, \beta} \int _{\partial \E} (\partial _\n \qbetau) \, \varphi_j^1 = \sum_{\beta =1}^{\alpha} \GS _{\alpha, \beta} \B_{\beta, j} = \GS(\alpha, 1:\alpha) \cdot \B(1:\alpha, j),\\
\end{split}
\]
where we used that it holds $\varphi_j^1|_{\partial \E} = \varphi_j^2|_{\partial \E}$.
\item Assume now $\varphi_j^2$ is an internal basis function. This case is a bit more involved. We write:
\begin{equation} \label{matrix Bbar internal element}
(\nabla \qalphad, \nabla \varphi_j^2)_{0,\E} = - \int _\E (\Delta \qalphad) \varphi _j^2.
\end{equation}
We expand $\Delta \qalphad$ into a combination of elements of the basis $\{\qbetad\}_{\beta=0}^{\npmd}$, since the Laplace operator eliminates the high ($\p-1$ and $\p$) polynomial degree contributions. We get:
\begin{equation} \label{expansion of the laplacian}
\Delta \qalphad = \sum_{\vert \beta \vert =0}^{\p-2} \Fbar_{\alpha,\beta} \qbetad.
\end{equation}
We only need to compute the entries of matrix $\Fbar$. For the purpose, we test \eqref{expansion of the laplacian} with $\qgammad$ thus obtaining:
\[
\Fbar_{\alpha,\gamma} = (\Delta \qalphad , \qgammad) _{0,\E} = - (\nabla \qalphad, \nabla \qgammad)_{0,\E} + \left(\partial_\n \qalphad, \qgammad \right)_{0,\partial \E}.
\]
The first term is nothing but $ - \Gtildebar_{\alpha,\gamma}$. We wonder how to compute the second term. We note that matrix $\L$ defined as:
\[
\L _{\alpha, \beta} = \int_{\partial \E} (\partial _\n\qalphau )\,\qbetau \quad \quad \forall \boldalpha,\, \boldbeta \in \mathbb N^2,\quad \vert \boldalpha\vert, \vert \boldbeta \vert = 0,\dots, \p,
\]
can be computed exactly. For the sake of completeness, we explicitly write how. Given $\mathcal E(\E)$ the set of edges of $\E$:
\[
\L _{\alpha, \beta} = \sum_{\e \in \mathcal E(\E)} \int_\e (\partial _\n\qalphau) \, \qbetau = \sum_{\e \in \mathcal E(\E)} \left\{ \sum_{k=0}^{\p}\omega_k^e \left( (\partial_\n \qalphau )\, \qbetau \right)(\nu_k^\e)   \right\},
\]
where $\omega_k^e$ and $\nu_k^e$, $k=0,\dots, \p$, are the $k$-th weights and nodes of the Gau\ss-Lobatto quadrature over edge $\e$. It is easy to check that if we set:
\[
\Lbar_{\alpha,\beta} = \int_{\partial \E} (\partial _\n \qalphad) \, \qbetad,
\]
then:
\[
\Lbar= \GS \cdot \L \cdot \GS^T.
\]
As a consequence:
\[
\Fbar = \Lbar - \Gtildebar.
\]
Now, we plug \eqref{expansion of the laplacian} in \eqref{matrix Bbar internal element} obtaining:
\[
\Bbar _{\alpha,j} = (\nabla \qalphad, \nabla \varphi_j^2) = - \sum_ {\beta=1}^{\npmd} \Fbar_{\alpha, \beta} (\qbetad, \varphi _j^2)_{0,\E} = - \sum _{\beta=1}^{\npmd} \Fbar_{\alpha,\beta} \Cbar _{\beta,j} = 
					-\Fbar(\alpha, 1:\npmd) \cdot \Cbar(1:\npmd,j),
\] 
where $\Cbar$ is a matrix defined as follows:
\[
\Cbar_{\alpha, j} = (\qalphad, \varphi^2_j)_{0,\E} \quad \quad \forall \boldalpha \in \mathbb N^2,\quad \vert \boldalpha \vert = 0,\dots, \p-2,\quad \forall j=1,\dots \NdofE.
\]
One has:
\[
(\qalphad, \varphi^2_j)_{0,\E} =
\begin{cases}
0 & \text{if } \varphi_j^2 \text{ is a boundary basis element}\\
\vert \E \vert \frac{1}{\vert \E \vert} \int_\E \qalphad \, \varphi_j ^2 = \delta_{\alpha, j} \vert \E \vert & \text{otherwise}
\end{cases}.
\]
\end{itemize}

{\footnotesize
\bibliography{bibliogr}
}
\bibliographystyle{ieeetr}
\bibliographystyle{siam}

\end{document}